\documentclass[12pt,a4paper,onecolumn]{article}
\usepackage{textcomp}
\usepackage{tikz}
\usepackage{amssymb}
\usepackage{caption}
\usepackage{color}
\usepackage{verbatim}
\usepackage{fancybox}
\usepackage{titlesec} 
\usepackage{mathtools}
\usepackage{amsmath}
\usepackage{esvect}
\usepackage{cite}
\usepackage{setspace}
\usepackage{slashbox}
\usepackage{booktabs}
\usepackage{diagbox}
\numberwithin{equation}{section}
\usepackage[hidelinks]{hyperref}
\usepackage{mathrsfs}
\usetikzlibrary{positioning,fit,calc}
\usepackage{array}
\usepackage{bm}
\usepackage{multirow}
\usepackage{hhline}
\usepackage[super]{nth}
\usetikzlibrary{positioning,decorations.pathreplacing}

\usepackage{mathrsfs}
\usetikzlibrary{positioning,fit,calc}
\usepackage[a4paper,left=1in,top=1in,right=1in,bottom=1in,nohead]{geometry} 
\usepackage[hidelinks]{hyperref}
\usepackage{array}
\usepackage{bm}
\usepackage[symbol]{footmisc} 
\numberwithin{equation}{section}
\def\correspondingauthor{\footnote{Corresponding author. Email: williewong088@gmail.com.}}

\tikzset{block/.style={draw,thick,text width=2cm,minimum height=1cm,align=center},
         line/.style={-latex}}
\newcolumntype{P}[1]{>{\centering\arraybackslash}m{#1}} 

\titleformat{\section}[block]{\large\scshape\bfseries}{\thesection.}{1em}{} 
\titleformat{\subsection}[block]{\bfseries}{\thesubsection.}{1em}{} 

\newtheorem{defn}{Definition}[section]
\newtheorem{thm}[defn]{Theorem}

\newtheorem{ppn}[defn]{Proposition}

\newtheorem{cor}[defn]{Corollary}
\newtheorem{lem}[defn]{Lemma}
\newtheorem{con}[defn]{Conjecture}

\begin{document}
\pagenumbering{arabic}
\begin{center}
    \textbf{\Large Kruskal-Katona's function and\\a variation of cross-intersecting antichains}
\vspace{0.1 in} 
    \\{\large W.H.W. Wong\correspondingauthor{}, E.G. Tay}
\vspace{0.1 in} 
\\National Institute of Education\\Nanyang Technological University, Singapore
\end{center}

\begin{abstract}
We prove some properties of the Kruskal-Katona function, and apply to the following variation of cross-intersecting antichains. Let $n\ge 4$ be an even integer and $\mathscr{A}$ and $\mathscr{B}$ be two cross-intersecting antichains of $\mathbb{N}_n$ with at most $k$ disjoint pairs, i.e. for all $A_i\in \mathscr{A}$, $B_j\in\mathscr{B}$, $A_i\cap B_j=\emptyset$ only if $i=j\le k$. We prove a best possible upper bound on $|\mathscr{A}|+|\mathscr{B}|$. Furthermore, we show that the extremal families contain only $\frac{n}{2}$ and $(\frac{n}{2}+1)$-sets.
\end{abstract}
\section{Introduction}

\indent\par Let $\mathbb{N}_n:=\{1,2,\ldots,n\}$ and $2^{\mathbb{N}_n}$ denote its power set for $n\in \mathbb{Z}^+$. For any integer $k$, $0\le k\le n$, ${{\mathbb{N}_n}\choose{k}}$ denotes the collection of all $k$-sets of $\mathbb{N}_n$. For any family $\mathscr{A}\subseteq 2^{\mathbb{N}_n}$, $\mathscr{A}^{(k)}$ denotes the collection of $k$-sets in $\mathscr{A}$, i.e. $\mathscr{A}^{(k)}=\mathscr{A}\cap{{\mathbb{N}_n}\choose{k}}$. $\mathscr{A}$ is said to be \textit{intersecting} if $X\cap Y\neq \emptyset$ for all $X,Y\in \mathscr{A}$. Two families $\mathscr{A},\mathscr{B}\subseteq 2^{\mathbb{N}_n}$ are said to be \textit{cross $t$-intersecting} if $|A\cap B|\ge t$ for all $A\in\mathscr{A}$ and all $B\in\mathscr{B}$. If $t=1$, we simply say that $\mathscr{A}$ and $\mathscr{B}$ are cross-intersecting.
\indent\par A notion closely related to the intersection of sets is the containment of sets. Two subsets $X$ and $Y$ of $\mathbb{N}_n$ are said to be \textit{independent} if $X\not\subseteq Y$ and $Y\not\subseteq X$. If $X$ and $Y$ are independent, we may say that $X$ is independent of $Y$. An \textit{antichain} or \textit{Sperner family} $\mathscr{A}$ of $\mathbb{N}_n$ is a collection of pairwise independent subsets of $\mathbb{N}_n$, i.e. for all $X,Y\in \mathscr{A}$, $X\not\subseteq Y$.
\indent\par The Erd{\" o}s-Ko-Rado's Theorem is central to the study of intersecting family of sets. It was later extended by Hilton and Milner to cross-intersecting families of sets.

\begin{thm} (Erd{\"o}s, Ko and Rado \cite{EKR})
Let $n\in\mathbb{Z}^+$ and $\mathscr{A}\subseteq {{\mathbb{N}_n}\choose{k}}$ be an intersecting family for some integer $k\le \frac{n}{2}$. Then, $|\mathscr{A}|\le {{n-1}\choose{k-1}}$.
\end{thm}

\begin{thm} (Hilton and Milner \cite{HM})
Let $n\in\mathbb{Z}^+$ and $\mathscr{A}, \mathscr{B}\subseteq {{\mathbb{N}_n}\choose{k}}$ be nonempty cross-intersecting families for some integer $k\le \frac{n}{2}$. Then, $|\mathscr{A}|+|\mathscr{B}|\le {{n}\choose{k}}-{{n-k}\choose{k}}+1$.
\end{thm}
\indent\par The above results saw extensions in various forms; Milner obtained what has now become a well-known analogue for an intersecting antichain. Frankl and Wong \cite{FP WHW} and Ou \cite{OY} independently obtained analogues on cross $t$-intersecting antichains.
\begin{thm}(Milner \cite{MEC})
Let $n,k\in\mathbb{Z}^+$ and $\mathscr{A}$ be an intersecting antichain of $\mathbb{N}_n$, where $|X\cap Y|\ge k$ for all $X,Y\in\mathscr{A}$. Then, $|\mathscr{A}|\le {{n}\choose{\lfloor \frac{n+k+1}{2}\rfloor}}$.
\end{thm}
\begin{thm}(Frankl and Wong \cite{FP WHW}, Ou \cite{OY})
Let $\mathscr{A}$ and $\mathscr{B}$ be two cross $t$-intersecting antichains of $\mathbb{N}_n$. Then,
\begin{align*}
|\mathscr{A}|+|\mathscr{B}|\le \max\limits_{t\le i\le \frac{n+t}{2}}\Bigg\{{{n}\choose{i}}+{{n}\choose{n+t-i}}\Bigg\}.
\end{align*}
with equality if and only if $\{\mathscr{A},\mathscr{B}\}=\{{{X}\choose{i^*}}, {{X}\choose{n+t-i^*}}\}$ for some integer $t\le i^*\le \frac{n+t}{2}$.
\end{thm}
\indent\par In this paper, we consider antichains of $\mathbb{N}_n$, $\mathscr{A}$ and $\mathscr{B}$, with at most $k$ disjoint pairs. Before introducing our results formally, we need to mention some classical results. 
\begin{thm}(Sperner \cite{SE}) For any $n\in \mathbb{Z}^+$, if $\mathscr{A}$ is an antichain of $\mathbb{N}_n$, then $|\mathscr{A}|\le {{n}\choose{\lfloor{n/2}\rfloor}}$. Furthermore, equality holds if and only if all members in $\mathscr{A}$ have the same size, ${\lfloor\frac{n}{2}\rfloor}$ or ${\lceil{\frac{n}{2}}\rceil}$.
\end{thm}
\indent\par For a family $\mathscr{A}\subseteq {{\mathbb{N}_n}\choose{k}}$, the \textit{shadow} and \textit{shade} of $\mathscr{A}$ are defined as
\begin{align*}
&\Delta \mathscr{A}:=\{ X\subseteq \mathbb{N}_n|\ |X|=k-1, X\subset Y\text{ for some } Y\in \mathscr{A}\},\text{ if }k>0,\text{ and }\\
&\nabla \mathscr{A}:=\{ X\subseteq \mathbb{N}_n|\ |X|=k+1, Y\subset X\text{ for some }Y\in \mathscr{A}\},\text{ if }k<n 
\end{align*}
respectively.
\indent\par We recall the following elementary inequalities due to Sperner,
\begin{lem} \cite{SE}\label{lemA8.1.6}
Let $\mathscr{A}$ be a collection of $k$-sets of $\mathbb{N}_n$. Then, 
\begin{align}
&|\nabla \mathscr{A}|/|\mathscr{A}|\ge \frac{n-k}{k+1},\label{eqA5.1.1} \text{ if }k<n, \text{ and }\\
&|\Delta \mathscr{A}|/|\mathscr{A}|\ge \frac{k}{n-k+1},\label{eqA5.1.2} \text{ if }k>0.
\end{align}
Furthermore, equality holds if and only if $\mathscr{A}=\emptyset$ or $\mathscr{A}={{\mathbb{N}_n}\choose{k}}$.
\end{lem}
\indent\par Let $\mathscr{F}\subseteq \mathbb{N}_n$. Define the top and bottom sizes $t(\mathscr{F})=\max\{|F|: F\in\mathscr{F}\}$ and $b(\mathscr{F})=\min\{|F|: F\in\mathscr{F}\}$. Following Sperner, let us define two new families obtained from $\mathscr{A}$.
\begin{align*}
\mathscr{F}_\circ=\mathscr{F}\backslash\mathscr{F}^{(t(\mathscr{F}))}\cup\Delta(\mathscr{F}^{(t(\mathscr{F}))}),\\
\mathscr{F}^\circ=\mathscr{F}\backslash\mathscr{F}^{(b(\mathscr{F}))}\cup\nabla(\mathscr{F}^{(b(\mathscr{F}))}).
\end{align*}
\indent\par The next lemma is well-known.
\begin{lem}\label{lemA5.1.7}
If $\mathscr{F}$ is a nonempty antichain $(\mathscr{F}\neq\{\emptyset\}$ or $\mathbb{N}_n)$, then both $\mathscr{F}_\circ$ and $\mathscr{F}^\circ$ are antichains and $t(\mathscr{F}_\circ)=t(\mathscr{F})-1$, $b(\mathscr{F}^\circ)=b(\mathscr{F})+1$.
\end{lem}
\indent\par The above bounds for the shadow $\Delta \mathscr{A}$ and shade $\nabla \mathscr{A}$ are not tight except for $\mathscr{A}=\emptyset$ or $\mathscr{A}={{\mathbb{N}_n}\choose{k}}$. A tight lower bound is given by the celebrated Kruskal-Katona's Theorem (KKT). KKT is closely related to the \textit{squashed order} of the $k$-sets. The squash relations $\le_s$ and $<_s$ are defined as follows. For $A,B\in {{\mathbb{N}_n}\choose{k}}$, $A\le_s B$ if the largest element of the symmetric difference $A+B:=(A-B)\cup (B-A)$ is in $B$. Furthermore, denote $A <_s B$ if $A\le_s B$ and $A\neq B$. For e.g., the $3$-subsets of $\mathbb{N}_5$ in squashed order are: $\bm{123} <_s \bm{124} <_s \bm{134} <_s \bm{234} <_s \bm{125} <_s \bm{135} <_s \bm{235} <_s \bm{145} <_s \bm{245} <_s \bm{345}$. Here, we omit the braces and write $\bm{abc}$ to represent the set $\{abc\}$, if there is no ambiguity. It can be easily shown that $<_s$ is anti-symmetric and transitive.
\indent\par We shall denote the collections of the first $m$ and last $m$ $k$-subsets of $\mathbb{N}_n$ in squashed order by $F_{n,k}(m)$ and $L_{n,k}(m)$ respectively. We use $C_{n,k}(m)$ to denote some collection of consecutive $k$-subsets of $\mathbb{N}_n$ in squashed order. We denote by $N^r_{n,k}(m)$, the collection $C_{n,k}(m)$ that follows $F_{n,k}(r)$ in squashed order. Then, KKT says that the shadow of a family $\mathscr{A}$ of $k$-sets is at least the size of the shadow of the first $|\mathscr{A}|$ $k$-sets in squashed order.
\begin{thm}(Kruskal \cite{KJB}, Katona \cite{KGOH}, and Clements and Lindstr{\"o}m \cite{CGF LB})
~\\Let $\mathscr{A}$ be a collection of $k$-sets of $\mathbb{N}_n$ and suppose the $k$-binomial representation of $|\mathscr{A}|$ is
\begin{align*}
|\mathscr{A}|={{a_k}\choose{k}}+{{a_{k-1}}\choose{k-1}}+\ldots+{{a_t}\choose{t}},
\end{align*}
where $a_k>a_{k-1}>\ldots>a_t\ge t\ge 1$. Then,
\begin{align*}
&|\Delta \mathscr{A}|\ge|\Delta F_{n,k}(|\mathscr{A}|)|={{a_k}\choose{k-1}}+{{a_{k-1}}\choose{k-2}}+\ldots+{{a_t}\choose{t-1}}.
\end{align*}
\end{thm}
\indent\par Lieby \cite{LP} proved that the shadow of the first $m$ $k$-sets of $\mathbb{N}_n$ in squashed order has the same cardinality as the shade of the last $m$ $(n-k)$-subsets of $\mathbb{N}_n$ in squashed order.
\begin{lem} (Lieby \cite{LP}) \label{lemA5.1.9}
For any integer $0\le m\le {{n}\choose{k}}$, $|\Delta F_{n,k}(m)|=|\nabla L_{n,n-k}(m)|$.
\end{lem}
\indent\par Let $S$ be a $k$-subset of $\mathbb{N}_n$. The \textit{new-shadow} and \textit{new-shade} of $S$ are defined as $\Delta_N S:=\{ X|\ X\in \Delta S, X\not\in \Delta T \text{ for all }T<_s S\}$ and $\nabla_N S:=\{ X |\ X\in \nabla S, X\not\in \nabla T \text{ for all }T>_s S\}$. Furthermore, if $\mathscr{A}$ is a collection of $k$-sets of $\mathbb{N}_n$, then the \textit{new-shadow} and \textit{new-shade} of $\mathscr{A}$ are defined as $\Delta_N \mathscr{A}:=\bigcup\limits_{S\in \mathscr{A}} \Delta_N S$ and $\nabla_N \mathscr{A}:=\bigcup\limits_{S\in \mathscr{A}} \nabla_N S$ respectively. We refer the interested readers to \cite{AI} and \cite{LP} for more details.
\indent\par It was shown by Clements \cite{CGF} that the size of the new shadow of any consecutive $m$ $k$-sets is at least that of the new shadow of the last $m$ $k$-sets in squashed order and its corresponding dual follows.
\begin{thm} (Clements \cite{CGF})\label{lemA5.1.10}
For any integer $0\le m\le {{n}\choose{k}}$, 
\begin{align*}
|\Delta_N C_{n,k}(m)|\ge |\Delta_N L_{n,k}(m)| \text{ and }|\nabla_N C_{n,k}(m)|\ge |\nabla_N F_{n,k}(m)|.
\end{align*}
\end{thm}
\section{Main results and motivation} 
\indent\par Acknowledging the contribution of KKT, we define the following functions $\kappa_{n,r}(\cdot)$ and $\kappa_{n,r}^*(\cdot)$.
\begin{defn}
Let $n,r$ be integers. Define
\begin{align*}
&\kappa_{n,r}(i):=|\Delta F_{n,r}(i)|-i \text{ and }\kappa_{n,r}^*(i):=\min\limits_{0\le j\le i}\kappa_{n,r}(j).
\end{align*}
\end{defn}
\indent\par This definition is due to $F_{n,r}(i)$ having the smallest possible shade, as shown by KKT. $\kappa_{n,r}(\cdot)$ is also known as the Kruskal-Katona function (KKF), which was studied by Frankl et al. \cite{FP MN RIZ TN}, and Minabutdinov and Manaev \cite{MAR MIE}. For any integer $n$, $\kappa_{n,r}(i)$ and $\kappa_{n.r}^*(i)$ can be computed using the $r$-binomial representation of $i$ and KKT. However, significant complexity resides in its computations due to its dependence on binomial representations. Despite so, we derive some useful properties of the KKF. We shall omit the subscript $n$ when there is no ambiguity.
\begin{ppn}\label{ppnA5.2.2}
$\kappa_{r}(m)<0$ if and only if $m\ge 1+\sum\limits_{i=1}^{r}{{2i-1}\choose{i}}$.
\end{ppn}
\begin{thm}\label{thmA5.2.3}
$\kappa_{r}(m)=\kappa_{r}^*(m)$ if and only if $m=\sum\limits_{i=t}^{r}{{a_i}\choose{i}}$ with $a_i\ge 2i-1$ for all $i=t,t+1,\ldots, r$.
\end{thm}
\begin{ppn}\label{ppnA5.2.4}
Let $n$, $a$ and $k$ be integers such that $0\le a,k \le {{n}\choose{\lceil n/2\rceil}}$. Then, 
\begin{align*}
\kappa_{\lceil \frac{n}{2}\rceil}({{n}\choose{\lceil n/2\rceil}})+\kappa_{\lceil \frac{n}{2}\rceil}^*(k)\le \kappa_{\lceil \frac{n}{2}\rceil}(a)+\kappa_{\lceil \frac{n}{2}\rceil}^*(k+{{n}\choose{\lceil n/2\rceil}}-a).
\end{align*}
\end{ppn}
\indent\par We apply the above results to the following variation of cross-intersecting antichains. Let $\mathscr{A}$ and $\mathscr{B}$ be antichains of $\mathbb{N}_n$ with at most $k$ disjoint pairs, i.e., for all $A_i\in \mathscr{A}$, $B_j\in\mathscr{B}$, $A_i\cap B_j=\emptyset$ only if $i=j\le k$. What is the maximum possible $|\mathscr{A}|+|\mathscr{B}|$?
\indent\par If $n$ is odd, it trivially follows from Sperner's Theorem that $|\mathscr{A}|+|\mathscr{B}|\le 2{{n}\choose{\lceil n/2\rceil}}$, with equality if $\mathscr{A}=\mathscr{B}={{\mathbb{N}_n}\choose{\lceil n/2\rceil}}$. Hence, we consider even integers $n$ and determine the following sharp bound.

\begin{thm}\label{thmA5.2.5}
Let $n\ge 4$ be an even integer and $\mathscr{A}$ and $\mathscr{B}$ be two antichains of $\mathbb{N}_n$. Suppose for some integer $k\le \min\{|\mathscr{A}|, |\mathscr{B}|\}$, and for all $A_i\in \mathscr{A}$, $B_j\in\mathscr{B}$, $A_i\cap B_j=\emptyset$ only if $i=j\le k$ (i.e. there are at most $k$ non-intersecting $A-B$ pairs). Then,  
\begin{align*}
|\mathscr{A}|+|\mathscr{B}|\le {{n}\choose{n/2}}+{{n}\choose{(n/2)+1}}-\kappa_{\frac{n}{2}}^*(k),
\end{align*}
where $\kappa_{\frac{n}{2}}^*(k)=0$ if $k<1+\sum\limits_{i=1}^{n/2}{{2i-1}\choose{i}}$ and $\kappa_{\frac{n}{2}}^*(k)<0$ otherwise.
Furthermore, equality holds if 
\\(i) $k<1+\sum\limits_{i=1}^{n/2}{{2i-1}\choose{i}}$, $\mathscr{A}={{\mathbb{N}_n}\choose{n/2}}$ and $\mathscr{B}={{\mathbb{N}_n}\choose{(n/2)+1}}$, or 
\\(ii) $k\ge 1+\sum\limits_{i=1}^{n/2}{{2i-1}\choose{i}}$, $\mathscr{A}={{\mathbb{N}_n}\choose{n/2}}$ and $\mathscr{B}=L_{n,\frac{n}{2}}(m)\cup {{\mathbb{N}_n}\choose{(n/2)+1}}-\nabla L_{n,\frac{n}{2}}(m)$, where $m\le k$ is an integer such that $\kappa_{\frac{n}{2}}^*(k)=\kappa_{\frac{n}{2}}(m)$.
\end{thm}
\indent\par We further prove that all pairs of extremal families contain only $\frac{n}{2}$ and $(\frac{n}{2}+1)$-sets.
\begin{thm}\label{thmA5.2.6}
Let $\mathscr{A}, \mathscr{B}$ and $k$ be as given in Theorem \ref{thmA5.2.5}. If $|\mathscr{A}|+|\mathscr{B}|={{n}\choose{n/2}}+{{n}\choose{(n/2)+1}}-\kappa_{\frac{n}{2}}^*(k)$, then for $\mathscr{X}=\mathscr{A},\mathscr{B}$, 
\\(i) $\mathscr{X}\subseteq {{\mathbb{N}_n}\choose{n/2}}\cup {{\mathbb{N}_n}\choose{(n/2)+1}}$,
\\(ii) $\mathscr{X}=\mathscr{X}^{(n/2)}\cup {{\mathbb{N}_n}\choose{(n/2)+1}}-\nabla\mathscr{X}^{(n/2)}$, and
\\(iii) $|\nabla\mathscr{X}^{(n/2)}|=|\nabla L_{n,\frac{n}{2}}(|\mathscr{X}^{(n/2)}|)|$.
\end{thm}
\indent\par Closely akin to this variation is the one with \textit{exactly} $k$ disjoint pairs, for which Theorem \ref{thmA5.2.5} also provides an upper bound. Moreover, it is tight if $\kappa_{\frac{n}{2}}(k)=\kappa_{\frac{n}{2}}^*(k)$ since $\mathscr{A}={{\mathbb{N}_n}\choose{n/2}}$ and $\mathscr{B}=L_{n,\frac{n}{2}}(k)\cup {{\mathbb{N}_n}\choose{(n/2)+1}}-\nabla L_{n,\frac{n}{2}}(k)$ make a pair of extremal antichains. Note that Theorem \ref{thmA5.2.3} gives the values of $k$ satisfying $\kappa_{\frac{n}{2}}(k)=\kappa_{\frac{n}{2}}^*(k)$. Hence, it is easy to see Corollary \ref{corA5.2.7} follow from the last two theorems.
\begin{cor}\label{corA5.2.7}
Let $\mathscr{A}, \mathscr{B}$ and $k$ be as given in Theorem \ref{thmA5.2.5}, except now $A_i\cap B_j=\emptyset$ if and only if $i=j\le k$ now. If $\kappa_{\frac{n}{2}}^*(k)=\kappa_{\frac{n}{2}}(k)$, then
\begin{align}
|\mathscr{A}|+|\mathscr{B}|\le {{n}\choose{n/2}}+{{n}\choose{(n/2)+1}}-\kappa_{\frac{n}{2}}(k). \label{eqA5.2.1}
\end{align}
Furthermore, if equality holds, then for $\mathscr{X}=\mathscr{A},\mathscr{B}$,
\\(i) $\mathscr{X}\subseteq {{\mathbb{N}_n}\choose{n/2}}\cup {{\mathbb{N}_n}\choose{(n/2)+1}}$,
\\(ii) $\mathscr{X}=\mathscr{X}^{(n/2)}\cup {{\mathbb{N}_n}\choose{(n/2)+1}}-\nabla\mathscr{X}^{(n/2)}$, and
\\(iii) $|\nabla\mathscr{X}^{(n/2)}|=|\nabla L_{n,\frac{n}{2}}(|\mathscr{X}^{(n/2)}|)|$.
\end{cor}
\indent\par The primary motivation of investigating these variations is their intricate connection with optimal orientations of a special family of graphs, known as the $G$ vertex-multiplications. In 2000, Koh and Tay \cite{KKM TEG 8} introduced $G$ vertex-multiplications and extended the results on complete $n$-partite graphs. Koh and Tay \cite{KKM TEG 11} further studied tree vertex-multiplications and Ng and Koh \cite{NKL KKM} investigated cycle vertex-multiplications. Moreover, Wong and Tay \cite{WHW TEG 6, WHW TEG 7} use the main results in this paper to derive conditions for vertex-multiplications in $\mathscr{C}_0$ and $\mathscr{C}_1$ for all trees of diameter $4$.

\section{Properties of the Kruskal-Katona function}
\indent\par In \cite{WHW TEG 4}, Wong and Tay defined the following notion $D(n,r)$ and proved some useful results; we list them as Lemmas \ref{lemA5.3.2}(a) and \ref{lemA5.3.4} and further prove some properties.
\begin{defn}(Wong and Tay \cite{WHW TEG 4}) 
\\For any positive integers $r$ and $n$, define
\begin{align*}
D(n,r):= \left\{
  \begin{array}{@{}ll@{}}
     {{n}\choose{r-1}}-{{n}\choose{r}}, & \text{if}\ r\le n, \\
    0, & \text{otherwise}. \\
  \end{array}\right.
\end{align*}
\end{defn}

\begin{lem}\label{lemA5.3.2} For any positive integers $r$,$m$ and $n$, 
\\(a) $D(n,r)\ \substack{>\\=\\<}\ 0 \iff r\ \substack{>\\=\\<}\ \frac{n+1}{2}$.
\\(b) $D(n-1,r-1)+D(n-1,r)=D(n,r)$ for $1\le r \le n$.
\\(c) $D(n+1,r)<D(n,r)$ for $n\ge 2r-1$.
\\(d) If $n\ge 2r-1$ and $n> m$, then $D(n,r)<D(m,r)$.
\\(e) If $r\ge 2$, then $D(m,1)\le D(1,1)$ and $D(m,r)\le D(2r-2,r)$.
\end{lem}
\textit{Proof}:
\\(a) This follows from $D(n,r)={{n}\choose{r-1}}-{{n}\choose{r}}={{n}\choose{r-1}}(1-\frac{n+1-r}{r})={{n}\choose{r-1}}(\frac{2r-1-n}{r})$.
\\(b) $D(n-1,r-1)+D(n-1,r)=[{{n-1}\choose{r-2}}-{{n-1}\choose{r-1}}]+[{{n-1}\choose{r-1}}-{{n-1}\choose{r}}]=[{{n-1}\choose{r-2}}+{{n-1}\choose{r-1}}]-[{{n-1}\choose{r-1}}+{{n-1}\choose{r}}]
={{n}\choose{r-1}}-{{n}\choose{r}}=D(n,r)$.
\\(c) $D(n+1,r)-D(n,r)=D(n,r-1)<0$ by (b) and (a).
\\(d) If $m< 2r-1$, then $D(m,r)>0\ge D(n,r)$. If $m\ge 2r-1$, then $D(n,r)<D(m,r)$ by (c).
\\(e) It is easy to see $D(m,1)\le 0=D(1,1)$. If $m\ge 2r-1$, then $D(m,r)\le 0<D(2r-2,r)$ by (a). If $m\le 2r-2$, then by (b) and (a) respectively, $D(m,r)-D(m-1,r)=D(m-1,r-1)\ge 0$.
\begin{flushright}
$\Box$
\end{flushright}

\indent\par An identity we will use often in our proofs is Chu Shih-Chieh's Identity (CSC), which is also known as ``Hockey Stick Identity''. (See \cite{KKM TEG 0} for more details.)
\begin{lem} (Chu Shih-Chieh's Identity) For any $r,k\in \mathbb{N}$, ${{r}\choose{0}}+{{r+1}\choose{1}}+\ldots+{{r+k}\choose{k}}={{r+k+1}\choose{k}}$.
\end{lem}
\indent\par The next lemma is an analogue of CSC in $D(n,r)$.
\begin{lem} (Wong and Tay \cite{WHW TEG 4})\label{lemA5.3.4} For any positive integer $j\ge 2$, $\sum\limits_{r=1}^j{D(j-2+r,r)}=1$.
\end{lem}
\indent\par We shall further derive properties of $D(n,r)$. Lemma \ref{lemA5.3.5} shows that the sum of the smallest (in absolute sense) negative term in $r$-th column, namely $D(2r,r)$, and each largest term, namely $D(2i-2,i)$, in the $i$-th column for $i=1,2,\ldots,r-1$ is negative.
\begin{lem}\label{lemA5.3.5}
For all positive integers $r$, $D(2r,r)+\sum\limits_{i=1}^{r-1}D(2i-2,i)< 0$.
\end{lem}
\textit{Proof}: We use induction on $r$. It is easy to verify for $r=1,2,3$.
\indent\par Assume $D(2r,r)+\sum\limits_{i=1}^{r-1}D(2i-2,i)<0$ for some positive integer $r\ge 3$. Then,
\begin{align*}
&\ D(2r+2,r+1)+\sum\limits_{i=1}^{r}D(2i-2,i)\\
=&\ D(2r+2,r+1)+D(2r-2,r)+\sum\limits_{i=1}^{r-1}D(2i-2,i)\\
<&\ D(2r+2,r+1)+D(2r-2,r)-D(2r,r) &\text{ (by induction hypothesis)}\\
=&\ {{2r-2}\choose{r-3}}+{{2r-2}\choose{r-4}}-2{{2r-2}\choose{r}}&\text{ (after rearranging)}\\
<&\ 0.
\end{align*}
\begin{flushright}
$\Box$
\end{flushright}
\indent\par The next well-known lemma will be found useful.
\begin{lem}\label{lemA5.3.6}
$f(x)={{x}\choose{y}}$ is an increasing function for $x\ge y$.
\end{lem}
\textit{Proof of Proposition \ref{ppnA5.2.2}}: Let the $r$-binomial representation of $m$ be $m=\sum\limits_{i=t}^{r}{{a_i}\choose{i}}$, where $a_{r}>a_{r-1}>\ldots>a_t\ge t\ge 1$. Also, denote $p=1+\sum\limits_{i=1}^{r}{{2i-1}\choose{i}}$.
\\($\Leftarrow$) Since $m\ge p$, we have $a_r\ge 2r-1$. Otherwise, 
\begin{align}
m=\sum\limits_{r=t}^{r}{{a_i}\choose{i}}\le\sum\limits_{i=0}^{r-t}{{2r-2-i}\choose{r-i}}\le\sum\limits_{i=0}^{r-1}{{2r-2-i}\choose{r-i}}={{2r-1}\choose{r}}-1<p,\label{eqA5.3.1}
\end{align}
where the first inequality is due to $a_{r-i}\le 2r-2-i$ for $i=0,1,\ldots,r-t$, and Lemma \ref{lemA5.3.6}, and the second equality due to CSC. However, this contradicts $m\ge p$.
\\
\\Case 1. $a_r=2r-1$.
\indent\par If $a_i=2i-1$ for all $i=t,t+1,\ldots, r$, then $\kappa_r(m)=\sum\limits_{i=t}^{r}D(2i-1,i)=0$ by Lemma \ref{lemA5.3.2}(a) (with equality if and only if $a_i=2i-1$ for all $i=t,t+1,\ldots,r$). Hence, assume there exists some integer $i$, $t\le i\le r-1$ such that $a_i\neq 2i-1$, and let $l$ be the largest such integer.
\\
\\Claim 1: $a_l\ge 2l$.
\indent\par Suppose $a_l<2l-1$. Then, 
\begin{align}
\sum\limits_{i=t}^{l}{{a_i}\choose{i}}\le\sum\limits_{i=0}^{l-t}{{2l-2-i}\choose{l-i}}\le \sum\limits_{i=0}^{l-1}{{2l-2-i}\choose{l-i}}={{2l-1}\choose{l}}-1,\label{eqA5.3.2}
\end{align}
where the first inequality is due to $a_{l-i}\le 2l-2-i$ for $i=0,1,\ldots,l-t$, and Lemma \ref{lemA5.3.6}, and the equality due to CSC. Consequently,
\begin{align*}
m-(p-1)=\sum\limits_{i=t}^{r}{{a_i}\choose{i}}-\sum\limits_{i=1}^{r}{{2i-1}\choose{i}}=\sum\limits_{i=t}^{l}{{a_i}\choose{i}}-\sum\limits_{i=1}^{l}{{2i-1}\choose{i}}\le-\sum\limits_{i=1}^{l-1}{{2i-1}\choose{i}},
\end{align*}
where the second equality follows from $a_i=2i-1$ for all $i=l+1,l+2,\ldots, r$, the inequality due to (\ref{eqA5.3.2}). This contradicts $m\ge p$. 
\indent\par Therefore,
\begin{align*}
\kappa_r(m)=\sum\limits_{i=t}^{l}D(a_i,i)\le D(2l,l)+\sum\limits_{i=t}^{l-1}D(2i-2,i)\le D(2l,l)+\sum\limits_{r=1}^{l-1}D(2i-2,i)<0, 
\end{align*}
where the first equality follows from $D(a_i,i)=D(2i-1,i)=0$ for $i=l,l+1,\ldots, r$, by Lemma \ref{lemA5.3.2}(a), the first inequality due to Lemma \ref{lemA5.3.2}(d) and (e), the second inequality due to $D(2i-2,i)>0$ for all $i=1,2,\ldots, t-1$, by Lemma \ref{lemA5.3.2}(a), and the last inequality by Lemma \ref{lemA5.3.5}. 
\\
\\Case 2. $a_r\ge 2r$.
\indent\par By a deduction similar to above, we have
\begin{align*}
\kappa_r(m)&=\sum\limits_{i=t}^{r}D(a_i,i)\le D(2r,r)+\sum\limits_{i=t}^{r-1}D(2i-2,i)\le D(2r,r)+\sum\limits_{r=1}^{r-1}D(2i-2,i)<0.
\end{align*}
\\($\Rightarrow$) We consider two cases of $a_r$.
\\Case 1. $a_r\ge2r-1$.
\indent\par If $a_i=2i-1$ for all $i=t,t+1,\ldots, r$, then $\kappa_r(m)=\sum\limits_{i=t}^{r}D(a_i,i)=\sum\limits_{i=t}^{r}D(2i-1,i)=0$ by Lemma \ref{lemA5.3.2}(a). This contradicts $\kappa_r(m)>0$. Hence, there exists some integer $i$, $t\le i\le r$ such that $a_i\neq 2i-1$, and let $l$ be the largest such integer.
\indent\par If $a_l\ge 2l$, then
\begin{align}
\sum\limits_{i=t}^{l}{{a_i}\choose{i}}\ge {{a_l}\choose{l}}\ge {{2l}\choose{l}}>\sum\limits_{i=0}^{l-1}{{2l-1-i}\choose{l-i}}\ge\sum\limits_{i=1}^{l}{{2i-1}\choose{i}},\label{eqA5.3.3}
\end{align}
where the second and last inequalities follow from Lemma \ref{lemA5.3.6} and the third inequality by CSC. So,
\begin{align*}
m-(p-1)=\sum\limits_{i=t}^{r}{{a_i}\choose{i}}-\sum\limits_{i=1}^{r}{{2i-1}\choose{i}}=\sum\limits_{i=t}^{l}{{a_i}\choose{i}}-\sum\limits_{i=1}^{l}{{2i-1}\choose{i}}>0,
\end{align*}
where the second equality follows from $a_i=2i-1$ for all $i=l+1,l+2,\ldots, r$, and the last inequality by (\ref{eqA5.3.3}). That is, $m\ge p$.
\indent\par Now, suppose $a_l\le 2l-2$. By Lemma \ref{lemA5.3.2}(a), $D(a_l,l)>0$. If $D(a_i,i)\ge 0$ for all $i=t,t+1,\ldots, l-1$, then 
\begin{align}
\sum\limits_{i=t}^{l}D(a_i,i)>0. \label{eqA5.3.4}
\end{align}
With $D(a_i,i)=0$ for $i=l,l+1,\ldots, r$, we have $\kappa_r(m)=\sum\limits_{i=t}^{l}D(a_i,i)+\sum\limits_{i=l+1}^{i}D(a_i,i)>0$, a contradiction. Hence, $D(a_i,i)<0$ for some $t\le i\le l-1$. Let $s$ be the smallest integer such that $D(a_i,i)>0$ for all $i=s, s+1, \ldots, l$.
\\
\\Claim 2. $a_{s}=2s-2$.
\indent\par Since $D(a_{s},s)>0$, it follows from Lemma \ref{lemA5.3.2}(a) that $a_{s}<2s-1$. Suppose $a_{s}\le 2s-3$. Then, $a_{s-1}\le a_{s}-1\le 2(s-1)-2$. By Lemma \ref{lemA5.3.2}(a), $D(a_{s-1},s-1)>0$, which contradicts the minimality of $s$. So, $a_{s}=2s-2$.
\indent\par Now,
\begin{align}
\sum\limits_{i=t}^{l}{D(a_i,i)}\ge\sum\limits_{i=t}^{s}{D(a_i,i)}\ge\sum\limits_{i=t}^{s}{D(s-2+i,i)}\ge\sum\limits_{i=1}^{s}{D(s-2+i,i)}=1. \label{eqA5.3.5}
\end{align}
The first inequality is due to $D(a_i,i)> 0$ for all $i=s+1,s+2,\ldots, l$ (equality holds if $s=l$). The second inequality is due to $a_i\le s-2+i$ for $i=t,t+1,\ldots, s$, and Lemma \ref{lemA5.3.2}(d). If $t=1$, the third inequality follows immediately. And, if $t>1$, the third inequality follows from $D(s-2+i,i)\le 0$ for $r=1,2,\ldots, t-1,$ by Lemma \ref{lemA5.3.2}(a). Invoking Lemma \ref{lemA5.3.4}, we obtain the last equality. It follows that
\begin{align*}
\kappa_r(m)=\sum\limits_{i=t}^{r}{D(a_i,i)}=\sum\limits_{i=t}^{l}{D(a_i,i)}\ge 1,
\end{align*}
a contradiction.
\\
\\Case 2. $a_r\le 2r-2$.
\indent\par By Lemma \ref{lemA5.3.2}(a), $D(a_r,r)>0$. If $D(a_i,i)\ge 0$ for all $i=t,t+1,\ldots, r-1$, then $\kappa_r(m)=\sum\limits_{i=t}^{r}D(a_i,i)>0$, a contradiction.
\indent\par So, there exists some integer $i$, $t\le i\le r-1$ such that $D(a_i,i)<0$. Let $s$ be the smallest integer such that $D(a_i,i)>0$ for all $i=s, s+1, \ldots, r$. As in Claim 2, it can be shown that $a_{s}=2s-2$. Similarly, we have $a_i\le s-2+i$ for all $i=t,t+1,\ldots, s$.
\indent\par So, 
\begin{align*}
\kappa_r(m)&=\sum\limits_{i=t}^{r}{D(a_i,i)}\ge D(2s-2,s)+\sum\limits_{i=t}^{s-1}{D(a_i,i)}\ge D(2s-2,s)+\sum\limits_{i=t}^{s-1}{D(s-2+i,i)}\\
&\ge\sum\limits_{i=1}^{s}{D(s-2+i,i)}=1.
\end{align*}
The first inequality is due to $D(a_i,i)>0$ for $i=s+1,s+2,\ldots, r$ (equality holds if $s=l$). The second inequality is due to $a_i\le s-2+i$ for $i=t,t+1,\ldots, s$, and Lemma \ref{lemA5.3.2}(d). If $t=1$, the third inequality follows immediately. And, if $t>1$, the third inequality follows from $D(s-2+i,i)<0$ for $i=1,2,\ldots, t-1,$ by Lemma \ref{lemA5.3.2}(a). Invoking Lemma \ref{lemA5.3.4} obtains the last equality. This contradicts $\kappa_r(m)<0$.
\begin{flushright}
$\Box$
\end{flushright}
\indent\par By being more careful in the previous proof, we have the following corollary.
\begin{cor}\label{corA5.3.7}
$\kappa_r(m)=0$ if and only if $m=\sum\limits_{i=t}^{r}{{2i-1}\choose{i}}$ for some integer $t\le r$. 
\end{cor}

\textit{Proof of Theorem \ref{thmA5.2.3}}:
~\\Case 1. $m\le \sum\limits_{i=1}^{r}{{2i-1}\choose{i}}$.
\indent\par Since $\kappa_r(0)=0$ and by Proposition \ref{ppnA5.2.2}, $\kappa_r^*(m)\ge 0$, we have $\kappa_r^*(m)=0$. The statement now follows from Corollary \ref{corA5.3.7}.
\\
\\Case 2. $m>\sum\limits_{i=1}^{r}{{2i-1}\choose{i}}$.
\\$(\Rightarrow)$ As shown in (\ref{eqA5.3.1}), $m>\sum\limits_{i=1}^{r}{{2i-1}\choose{i}}$ implies $a_r\ge 2r-1$. Suppose  $a_i\ge 2i-1$ for all $i=l+1,l+2,\ldots, r$ and $a_l\le 2l-2$ for some integer $l$, where $t\le l<r$. Now, proceed as we did in Case 1 of Proposition \ref{ppnA5.2.2} to conclude $\sum\limits_{i=t}^{l}D(a_i,i)>0$ (see (\ref{eqA5.3.4}) and (\ref{eqA5.3.5})). Hence, noting also that $\sum\limits_{r=l+1}^{r}{{a_i}\choose{i}}<m$,
\begin{align*}
\kappa_r(m)=\sum\limits_{i=t}^{l}D(a_i,i)+\sum\limits_{r=l+1}^{r}D(a_i,i)>\kappa_r(\sum\limits_{r=l+1}^{r}{{a_i}\choose{i}})\ge \kappa_r^*(m),
\end{align*}
which contradicts $\kappa_{r}(m)=\kappa_{r}^*(m)$.
\\
\\$(\Leftarrow)$ Let $m^*<m$ with the $r$-binomial representation be $m^*=\sum\limits_{i=t^*}^{r}{{b_i}\choose{i}}$. It suffices to show that 
\begin{align}
\kappa_r(m)\le\kappa_r(m^*). \label{eqA5.3.6}
\end{align}
\indent\par Let $q=\max\{t,t^*\}$. If $a_i=b_i$ for all $i=q, q+1,\ldots, r$, then $m^*<m$ implies $t^*>t$. So, 
\begin{align*}
\kappa_r(m)=\sum\limits_{i=t}^{r}D(a_i,i)\le \sum\limits_{i=t^*}^{r}D(b_i,i)=\kappa_r(m^*),
\end{align*}
since $D(a_i,i)\le 0$ for $i=t,t+1,\ldots, t^*-1$, by Lemma \ref{lemA5.3.2}(a). That is, (\ref{eqA5.3.6}) holds.
\indent\par Suppose there exists some integer $i=q,q+1,\ldots, r$ such that $a_i\neq b_i$ and let $l$ be the largest such integer.
\\
\\Claim. $a_l>b_l$.
\indent\par Suppose $a_l<b_l$. Then,
\begin{align}
\sum\limits_{i=t^*}^{l}{{b_i}\choose{i}}\ge {{b_l}\choose{l}}\ge {{a_l+1}\choose{l}}>\sum\limits_{i=0}^{l-1}{{a_l-i}\choose{l-i}}\ge\sum\limits_{i=1}^{l}{{a_i}\choose{i}}\ge \sum\limits_{i=t}^{l}{{a_i}\choose{i}},\label{eqA5.3.7}
\end{align}
where the second and second last inequalities follow from Lemma \ref{lemA5.3.6} and the third inequality by CSC. So,
\begin{align*}
m^*-m=\sum\limits_{i=t^*}^{r}{{b_i}\choose{i}}-\sum\limits_{i=t}^{r}{{a_i}\choose{i}}=\sum\limits_{i=t^*}^{l}{{b_i}\choose{i}}-\sum\limits_{i=t}^{l}{{a_i}\choose{i}}>0,
\end{align*}
where the second equality follows from $a_i=2i-1$ for all $i=l+1,l+2,\ldots, r$, and the last inequality by (\ref{eqA5.3.7}). This contradicts $m^*<m$. Hence, the claim follows.
\indent\par Now, we consider two cases of $b_l$.
\\Case 1. $b_l\ge 2l-1$.
\indent\par Then, $b_{l-i}\le b_l-i$ for $i=1,2,\ldots, l-t^*$. Consequently, for $i=1,2,\ldots, l$, $b_l-i\ge 2(l-i)-1$ and 
\begin{align}
D(b_l-i,l-i)\le 0 \label{eqA5.3.8}
\end{align}
by Lemma \ref{lemA5.3.2}(a). Furthermore, for all $i=0,1,2,\ldots, l-t^*$,
\begin{align}
D(b_{l-i},l-i)\ge D(b_l-i,l-i) \label{eqA5.3.9}
\end{align}
by Lemma \ref{lemA5.3.2}(d). So, 
\begin{align*}
\sum\limits_{i=0}^{l-t^*}D(b_{l-i},l-i)\ge \sum\limits_{i=0}^{l-t^*}D(b_l-i,l-i)\ge \sum\limits_{i=0}^{l}D(b_l-i,l-i)=D(b_l+1,l)+1,
\end{align*}
where the first and second inequalities follow from (\ref{eqA5.3.9}) and (\ref{eqA5.3.8}) respectively. The equality is due to CSC. Now, invoking Lemma \ref{lemA5.3.2}(d) on $b_l<a_l$ and $a_l\ge 2l-1$ gives 
\begin{align}
\sum\limits_{i=t^*}^{l}D(b_i,r)\ge D(b_l+1,1)+1> D(a_l,l)+1. \label{eqA5.3.10}
\end{align}
\indent\par So,
\begin{align*}
\kappa_r(m)-\kappa_r(m^*)=&\sum\limits_{i=t}^{r}D(a_i,i)-\sum\limits_{i=t^*}^{r}D(b_i,i)=\sum\limits_{i=t}^{r}D(a_i,i)-\sum\limits_{i=t^*}^{l}D(b_i,i)-\sum\limits_{r=l+1}^{r}D(b_i,i)\\
<&\sum\limits_{i=t}^{r}D(a_i,i)-D(a_l,l)-1-\sum\limits_{r=l+1}^{r}D(b_i,i)=\sum\limits_{i=t}^{l-1}D(a_i,i)-1\\
\le& -1.
\end{align*}
The first inequality follows from (\ref{eqA5.3.10}) and the last equality due to $a_i=b_i$ for $i=l+1,l+2,\ldots, r$. The last inequality is due to $D(a_i,i)\le 0$ for all $i=t,t+1,\ldots, l-1$, by Lemma \ref{lemA5.3.2}(b) since $a_i\ge 2i-1$ for $i=t,t+1,\ldots,r$.
\\
\\Case 2. $b_l<2l-1$.
\indent\par If $D(b_i,i)\ge 0$ for all $i=t^*,t^*+1,\ldots, l$, then 
\begin{align*}
\kappa_r(m)-\kappa_r(m^*)=\sum\limits_{i=t}^{r}D(a_i,i)-\sum\limits_{i=t^*}^{r}D(b_i,i)=\sum\limits_{i=t}^{l}D(a_i,i)-\sum\limits_{i=t^*}^{l}D(b_i,i)\le 0,
\end{align*}
where the second equality is due to $a_i=b_i$ for $i=l+1,l+2,\ldots, \frac{n}{2}$, and the inequality due to $D(a_i,i)\le 0$ by Lemma \ref{lemA5.3.2}(a). So, we have (\ref{eqA5.3.6}) as required.
\indent\par Hence, we assume $D(b_i,i)<0$ for some $i=t^*,t^*+1,\ldots,l-1$. Let $s$ be the smallest integer such that $D(b_i,i)>0$ for all $i=s,s+1,\ldots, l$. As in Claim 2 of Proposition \ref{ppnA5.2.2}, we can prove that $b_s=2s-2$.
\indent\par Note for all $i=1,2,\ldots, s$, we have $s-2+i\ge 2i-1$ implying $D(s-2+i,i)\le 0$ by Lemma \ref{lemA5.3.2}(a). Furthermore, with $b_i\le s-2+i$ for all $i=t^*,t^*+1,\ldots, s$, it follows by Lemma \ref{lemA5.3.2}(d) that $D(b_i,i)\ge D(s-2+i,i)$ for $i=1,2,\ldots, s$. Consequently,
\begin{align}
\sum\limits_{i=t^*}^{s}D(b_i,i)\ge \sum\limits_{i=t^*}^{s} D(s-2+i,i)\ge \sum\limits_{i=1}^{s} D(s-2+i,i)=1. \label{eqA5.3.11}
\end{align}
\indent\par So,
\begin{align*}
\kappa_r(m)-\kappa_r(m^*)=&\sum\limits_{i=t}^{r}D(a_i,i)-\sum\limits_{i=t^*}^{r}D(b_i,i)=\sum\limits_{i=t}^{r}D(a_i,i)-\sum\limits_{i=t^*}^{s}D(b_i,i)-\sum\limits_{r=s+1}^{r}D(b_i,i)\\
\le&\sum\limits_{i=t}^{r}D(a_i,i)-1-\sum\limits_{r=s+1}^{r}D(b_i,i)\le\sum\limits_{i=t}^{s}D(a_i,i)-1\\
\le&-1.
\end{align*}
The first inequality is due (\ref{eqA5.3.11}) and the second inequality due to $D(b_i,i)>0\ge D(a_i,i)$ for $i=s+1,s+2,\ldots, l$, and $D(b_i,i)=D(a_i,i)$ for $i=l+1,l+2,\ldots, r$. We remark that $\sum\limits_{i=t}^{s}D(a_i,i)=0$ if $s<t$. Hence, (\ref{eqA5.3.6}) follows as desired.
\begin{flushright}
$\Box$
\end{flushright}

\begin{lem}\label{lemA8.3.8}
Let $m$ be an integer, where $0\le m \le{{n}\choose{\lceil n/2\rceil}}$. Then, $\kappa_{\lceil\frac{n}{2}\rceil}(m)\ge \kappa_{\lceil\frac{n}{2}\rceil}({{n}\choose{\lceil n/2\rceil}})$. Furthermore, equality holds if and only if $m={{n}\choose{n/2}}$.
\end{lem}
\textit{Proof}: Suppose $n$ is odd. By Proposition \ref{ppnA5.2.2}, $\kappa_{\lceil\frac{n}{2}\rceil}(m)\ge 0=\kappa_{\lceil\frac{n}{2}\rceil}({{n}\choose{\lceil n/2\rceil}})$ for all $0\le m \le{{n}\choose{\lceil n/2\rceil}}<1+\sum\limits_{i=1}^{\lceil n/2\rceil}{{2i-1}\choose{i}}$.
\indent\par Suppose $n$ is even. Noting that $\kappa_{\frac{n}{2}}({{n}\choose{n/2}})={{n}\choose{(n/2)+1}}-{{n}\choose{n/2}}$,
\begin{align*}
\kappa_{\frac{n}{2}}(m)=|\Delta F_{n,\frac{n}{2}}(m)|-m\ge \Big[\frac{n}{n+2}-1\Big]\cdot m\ge \Big[\frac{n}{n+2}-1\Big] {{n}\choose{n/2}}=\kappa_{\frac{n}{2}}({{n}\choose{n/2}}).
\end{align*}
where the first inequality follows from (\ref{eqA5.1.2}).
\indent\par Trivially, if $m={{n}\choose{n/2}}$, then $\kappa_{\frac{n}{2}}(m)=\kappa_{\frac{n}{2}}({{n}\choose{n/2}})$. If $\kappa_{\frac{n}{2}}(m)=\kappa_{\frac{n}{2}}({{n}\choose{n/2}})$, then equality must hold throughout, 
particularly the first inequality. By (\ref{eqA5.1.2}), this implies $m={{n}\choose{n/2}}$ or $m=0$. If $m=0$, then 
$\kappa_{\frac{n}{2}}(m)-\kappa_{\frac{n}{2}}({{n}\choose{n/2}})={{n}\choose{n/2}}-{{n}\choose{(n/2)+1}}>0$. So, it remains that $m={{n}\choose{n/2}}$.
\begin{flushright}
$\Box$
\end{flushright}

\textit{Proof of Proposition \ref{ppnA5.2.4}}: In this proof, we denote $\kappa_{\lceil\frac{n}{2}\rceil}$ (and $\kappa_{\lceil\frac{n}{2}\rceil}^*$) with $\bm{\kappa}$ (and $\bm{\kappa^*}$ resp.) for brevity. We apply induction on $k$. Consider $k=0$.
\\
\\Case B1. $a\le \sum\limits_{i=1}^{\lceil n/2 \rceil}{{2i-1}\choose{i}}$ and ${{n}\choose{\lceil n/2\rceil}}-a\le \sum\limits_{i=1}^{\lceil n/2 \rceil}{{2i-1}\choose{i}}$.
\indent\par By Proposition \ref{ppnA5.2.2}, $\bm{\kappa}(a)\ge 0$ and $\bm{\kappa}^*({{n}\choose{\lceil n/2\rceil}}-a)=0=\bm{\kappa}^*(0)$. So, $\bm{\kappa}({{n}\choose{\lceil n/2\rceil}})+\bm{\kappa}^*(0)-\bm{\kappa}(a)-\bm{\kappa}^*({{n}\choose{\lceil n/2\rceil}}-a)\le \bm{\kappa}({{n}\choose{\lceil n/2\rceil}})\le 0$.

\indent\par We remark that Cases B2 and B3 do not apply to odd integers $n$ since ${{n}\choose{\lceil n/2\rceil}}<\sum\limits_{i=1}^{\lceil n/2\rceil}{{2i-1}\choose{i}}$.
\\
\\Case B2. $a>\sum\limits_{i=1}^{n/2}{{2i-1}\choose{i}}$.
\indent\par Since $\sum\limits_{i=1}^{n/2}{{2i-1}\choose{i}}>\frac{1}{2}{{n}\choose{n/2}}$, it follows that ${{n}\choose{n/2}}-a<\frac{1}{2}{{n}\choose{n/2}}<\sum\limits_{i=1}^{n/2}{{2i-1}\choose{i}}$ and $\bm{\kappa}^*({{n}\choose{n/2}}-a)=0=\bm{\kappa^*}(0)$ by Proposition \ref{ppnA5.2.2}. Furthermore, $\bm{\kappa}({{n}\choose{n/2}})+\bm{\kappa}^*(0)-\bm{\kappa}(a)-\bm{\kappa}^*({{n}\choose{n/2}}-a)=\bm{\kappa}({{n}\choose{n/2}})-\bm{\kappa}(a)\le0$ by Lemma \ref{lemA8.3.8}.
\\
\\Case B3. ${{n}\choose{n/2}}-a>\sum\limits_{i=1}^{n/2}{{2i-1}\choose{i}}$.
\indent\par Since $\sum\limits_{i=1}^{n/2}{{2i-1}\choose{i}}>\frac{1}{2}{{n}\choose{n/2}}$, it follows that $a<\frac{1}{2}{{n}\choose{n/2}}<\sum\limits_{i=1}^{n/2}{{2i-1}\choose{i}}$ and $\bm{\kappa}(a)\ge 0$ by Proposition \ref{ppnA5.2.2}. Furthermore, $\bm{\kappa}({{n}\choose{n/2}})+\bm{\kappa}^*(0)-\bm{\kappa}(a)-\bm{\kappa}^*({{n}\choose{n/2}}-a)\le\bm{\kappa}({{n}\choose{n/2}})-\bm{\kappa}^*({{n}\choose{n/2}}-a)=\bm{\kappa}^*({{n}\choose{n/2}})-\bm{\kappa}^*({{n}\choose{n/2}}-a)\le 0$ by Lemma \ref{lemA8.3.8}.
\indent\par For the induction case, assume $\bm{\kappa}({{n}\choose{n/2}})+\bm{\kappa}^*(k)-\bm{\kappa}(a)-\bm{\kappa}^*(k+{{n}\choose{n/2}}-a)\le 0$ for some integer $k\ge 0$. We want to show $\bm{\kappa}({{n}\choose{n/2}})+\bm{\kappa}^*(k+1)-\bm{\kappa}(a)-\bm{\kappa}^*(k+1+{{n}\choose{n/2}}-a)\ge 0$.
\indent\par Let $m$ be the largest integer satisfying $m\le k+1$ and $\bm{\kappa}^*(k+1)=\bm{\kappa}(m)$, and let $b$ be an integer satisfying $b\le k+1+{{n}\choose{n/2}}-a$ and $\bm{\kappa}^*(k+1+{{n}\choose{n/2}}-a)=\bm{\kappa}(b)$.
\\
\\Case I1. $b< k+1+{{n}\choose{\lceil n/2\rceil}}-a$.
\indent\par Then, $\bm{\kappa}^*(k+1+{{n}\choose{\lceil n/2\rceil}}-a)=\bm{\kappa}(b)=\bm{\kappa}^*(k+{{n}\choose{\lceil n/2\rceil}}-a)$. Note that $\bm{\kappa}^*(k+1)\le \bm{\kappa}^*(k)$ by definition. It follows that

\begin{align*}
&\ \bm{\kappa}({{n}\choose{\lceil n/2\rceil}})+\bm{\kappa}^*(k+1)-\bm{\kappa}(a)-\bm{\kappa}^*(k+1+{{n}\choose{\lceil n/2\rceil}}-a)\\
\le&\ \bm{\kappa}({{n}\choose{\lceil n/2\rceil}})+\bm{\kappa}^*(k)-\bm{\kappa}(a)-\bm{\kappa}^*(k+{{n}\choose{\lceil n/2\rceil}}-a)\\
\le&\  0, 
\end{align*} by induction hypothesis.
\\
\\Case I2. $b= k+1+{{n}\choose{\lceil n/2\rceil}}-a$.
\indent\par Since $m\le k+1$, we have $b-m\ge {{n}\choose{\lceil n/2\rceil}}-a$. 
\\
\\Subcase I2.1. $m=k+1$, i.e. $b-m={{n}\choose{\lceil n/2\rceil}}-a$.
\indent\par Then,
\begin{align*}
&\ \bm{\kappa}({{n}\choose{\lceil n/2\rceil}})+\bm{\kappa}^*(k+1)-\bm{\kappa}(a)-\bm{\kappa}^*(k+1+{{n}\choose{\lceil n/2\rceil}}-a) \\
=&\ \bm{\kappa}({{n}\choose{\lceil n/2\rceil}})+\bm{\kappa}(m)-\bm{\kappa}(a)-\bm{\kappa}(b)\\
=&\ (b-m)-\Big[{{n}\choose{\lceil n/2\rceil}}-a\Big]+\Big[|\Delta F_{n,\frac{n}{2}}({{n}\choose{\lceil n/2\rceil}})|-|\Delta F_{n,\frac{n}{2}}(a)|\Big]-\Big[|\Delta F_{n,\frac{n}{2}}(b)|-|\Delta F_{n,\frac{n}{2}}(m)|\Big]\\
=&\ |\Delta_N L_{n,\frac{n}{2}}({{n}\choose{\lceil n/2\rceil}}-a)|-|\Delta_N N^m_{n,\frac{n}{2}}(b-m)|\\
=&\ |\Delta_N L_{n,\frac{n}{2}}(b-m)|-|\Delta_N N^m_{n,\frac{n}{2}}(b-m)|\\
\le&\ 0,
\end{align*}
where the last inequality follows from Theorem \ref{lemA5.1.10}.
\\
\\Subcase I2.2. $m<k+1$, i.e. $b-m>{{n}\choose{\lceil n/2\rceil}}-a$.
\indent\par Note that $m<k+1$ implies 
\begin{align}
|\Delta_N N^m_{n,\frac{n}{2}}(k+1-m)|> k+1-m. \label{eqA5.3.12}
\end{align}
For if $|\Delta_N N^m_{n,\frac{n}{2}}(k+1-m)|\le k+1-m$, then $\bm{\kappa}(k+1)=|\Delta F_{n,\frac{n}{2}}(k+1)|-(k+1)=|\Delta F_{n,\frac{n}{2}}(m)|-m+|\Delta_N N^m_{n,\frac{n}{2}}(k+1-m)|-(k+1-m)\le\bm{\kappa}(m)$, 
which contradicts the maximality of $m$. 
\indent\par So,
\begin{align*}
&\ \bm{\kappa}({{n}\choose{\lceil n/2\rceil}})+\bm{\kappa}^*(k+1)-\bm{\kappa}(a)-\bm{\kappa}^*(k+1+{{n}\choose{\lceil n/2\rceil}}-a)\\
=&\ \bm{\kappa}({{n}\choose{\lceil n/2\rceil}})+\bm{\kappa}(m)-\bm{\kappa}(a)-\bm{\kappa}(b)\\
=&\ -\Big[{{n}\choose{\lceil n/2\rceil}}-a-b\Big]-m+\Big[|\Delta F_{n,\frac{n}{2}}({{n}\choose{\lceil n/2\rceil}})|-|\Delta F_{n,\frac{n}{2}}(a)|\Big]-\Big[|\Delta F_{n,\frac{n}{2}}(b)|-|\Delta F_{n,\frac{n}{2}}(m)|\Big]\\
=&\ (k+1-m)+|\Delta_N L_{n,\frac{n}{2}}({{n}\choose{\lceil n/2\rceil}}-a)|-|\Delta_N N^{k+1}_{n,\frac{n}{2}}(b-(k+1))|-|\Delta_N N^m_{n,\frac{n}{2}}(k+1-m)|\\
<&\ 0,
\end{align*}
where the inequality is due to $b-(k+1)={{n}\choose{\lceil n/2\rceil}}-a$ and thus, $|\Delta_N N^{k+1}_{n,\frac{n}{2}}(b-(k+1))|\ge |\Delta_N L_{n,\frac{n}{2}}({{n}\choose{\lceil n/2\rceil}}-a)|$ by Theorem \ref{lemA5.1.10} and (\ref{eqA5.3.12}).
\begin{flushright}
$\Box$
\end{flushright}
\section{Proofs of Theorems \ref{thmA5.2.5} and \ref{thmA5.2.6}}
\indent\par To prove Theorem \ref{thmA5.2.5}, we shall employ Sperner's operations in a manner similar to that used by Sperner himself to prove Sperner's Theorem \cite{SE}.
\\\\\textit{Proof of Theorem \ref{thmA5.2.5}}: In view of (\ref{eqA5.1.1}), we may replace any element, $A_i$ ($B_i$ resp.), of size $<\frac{n}{2}$ in $\mathscr{A}$ ($\mathscr{B}$ resp.) with an equal number of $\frac{n}{2}$-sets, say $A^{\uparrow}_i$ ($B^{\uparrow}_i$ resp.), from $\nabla \mathscr{A}$ ($\nabla \mathscr{B}$ resp.). If $|A_i|=\frac{n}{2}$ ($|B_i|=\frac{n}{2}$ resp.), then $A^{\uparrow}_i:=A_i$ ($B^{\uparrow}_i:=B_i$) simply. Similarly, by (\ref{eqA5.1.2}), we may replace those $A_i$'s ($B_i$'s resp.) of size $>\frac{n}{2}+1$ by an equal number of $(\frac{n}{2}+1)$-sets, say $A^{\downarrow}_i$ ($B^{\downarrow}_i$ resp.) from their shadow. If $|A_i|=\frac{n}{2}+1$ ($|B_i|=\frac{n}{2}+1$ resp.), then $A^{\downarrow}_i:=A_i$ ($B^{\downarrow}_i:=B_i$) simply. By Lemma \ref{lemA5.1.7}, the replaced sets $\mathscr{A}$ and $\mathscr{B}$ are still antichains. Furthermore, $A_i\cap B_j \neq\emptyset\implies A^{\uparrow}_i\cap B^{\uparrow}_j \neq\emptyset$ since $A_i\subseteq A^{\uparrow}_i$ and $B_j\subseteq B^{\uparrow}_j$.
\indent\par Consequently, we may assume $\mathscr{X}\subseteq {{\mathbb{N}_n}\choose{n/2}}\cup {{\mathbb{N}_n}\choose{(n/2)+1}}$ for $\mathscr{X}=\mathscr{A},\mathscr{B}$. So, $|T|\le k$, where $T:=\{i|A^{\uparrow}_i\cap B^{\uparrow}_i=\emptyset\}$. Partition $\mathscr{X}$ into $\mathscr{X}_1:=X\cap {{\mathbb{N}_n}\choose{n/2}}$ and $\mathscr{X}_2:=\mathscr{X}-\mathscr{X}_1$. Then, $|\mathscr{A}_1|+|\mathscr{B}_1|\le |T|+{{n}\choose{n/2}} \le k+{{n}\choose{n/2}}$, which implies 
\begin{align}
|\mathscr{A}_1|+|\mathscr{B}_1|\le k+{{n}\choose{n/2}}.\label{eqA5.4.1}
\end{align}
\indent\par Now,
\begin{align}
&\ |\mathscr{A}|+|\mathscr{B}|\label{eqA5.4.2}\\
=&\ |\mathscr{A}_1|+|\mathscr{A}_2|+|\mathscr{B}_1|+|\mathscr{B}_2|\nonumber\\
\le&\ |\mathscr{A}_1|+|{{\mathbb{N}_n}\choose{(n/2)+1}}-\nabla \mathscr{A}_1|+|\mathscr{B}_1|+|{{\mathbb{N}_n}\choose{(n/2)+1}}-\nabla \mathscr{B}_1|\nonumber\\
\le&\ 2{{n}\choose{(n/2)+1}}+|\mathscr{A}_1|-|\nabla L_{n,\frac{n}{2}}(|\mathscr{A}_1|)|+|\mathscr{B}_1|-|\nabla L_{n,\frac{n}{2}}(|\mathscr{B}_1|)|\nonumber\\
=&\ 2{{n}\choose{(n/2)+1}}-\kappa(|\mathscr{A}_1|)-\kappa(|\mathscr{B}_1|)\nonumber\\
\le&\ 2{{n}\choose{(n/2)+1}}-\kappa(|\mathscr{A}_1|)-\kappa^*(k+{{n}\choose{n/2}}-|\mathscr{A}_1|)\nonumber\\
\le&\ 2{{n}\choose{(n/2)+1}}-\kappa({{n}\choose{(n/2)}})-\kappa^*(k)\nonumber\\
=&\ {{n}\choose{(n/2)}}+{{n}\choose{(n/2)+1}}-\kappa^*(k).\nonumber
\end{align}
where the second inequality is due to $|\nabla \mathscr{X}_1|\ge |\nabla L_{n,\frac{n}{2}}(|\mathscr{X}_1|)|=|\Delta F_{n,\frac{n}{2}}(|\mathscr{X}_1|)|$ for both $\mathscr{X}=\mathscr{A},\mathscr{B}$, by KKT and Lemma \ref{lemA5.1.9}. The fourth and fifth inequalities follow from (\ref{eqA5.4.1}) and Proposition \ref{ppnA5.2.4} respectively.
\begin{flushright}
$\Box$
\end{flushright}
\textit{Proof} of Theorem \ref{thmA5.2.6}: In this proof,  $\mathscr{X}$ and $\mathscr{X}_1\cup\mathscr{X}_2$ denote the antichain before and after Sperner operations respectively for $\mathscr{X}=\mathscr{A},\mathscr{B}$. Since $|\mathscr{A}|+|\mathscr{B}|={{n}\choose{n/2}}+{{n}\choose{(n/2)+1}}+\kappa^*(k)$, equality must hold through in (\ref{eqA5.4.2}). Particularly, we must have 
\begin{align}
&\kappa(|\mathscr{A}_1|)+\kappa^*(k+{{n}\choose{n/2}}-|\mathscr{A}_1|)=\kappa({{n}\choose{n/2}})+\kappa^*(k),  \label{eqA5.4.3}\\
&\mathscr{A}_2={{\mathbb{N}_n}\choose{(n/2)+1}}-\nabla\mathscr{A}_1, \text{ }
\mathscr{B}_2={{\mathbb{N}_n}\choose{(n/2)+1}}-\nabla\mathscr{B}_1, \label{eqA5.4.4}\\
&|\nabla \mathscr{A}_1|=|\nabla L_{n,\frac{n}{2}}(|\mathscr{A}_1|), \text{ and }
|\nabla \mathscr{B}_1|=|\nabla L_{n,\frac{n}{2}}(|\mathscr{B}_1|). \label{eqA5.4.5}
\end{align}
\indent\par From
\begin{align}
|\mathscr{A}|+|\mathscr{B}|={{n}\choose{n/2}}+{{n}\choose{(n/2)+1}}-\kappa^*(k)\ge {{n}\choose{n/2}}+{{n}\choose{(n/2)+1}}\label{eqA5.4.6}
\end{align}
and by Sperner's Theorem, $|\mathscr{A}|\le {{n}\choose{n/2}}$ and $|\mathscr{B}|\le {{n}\choose{n/2}}$, it follows that $|\mathscr{A}|>0$ and $|\mathscr{B}|>0$.
\indent\par We first show $t(\mathscr{X})\le\frac{n}{2}+1$ for $\mathscr{X}=\mathscr{A},\mathscr{B}$. Suppose $\mathscr{A}^{(s)}\neq\emptyset$ for some $s>\frac{n}{2}+1$. For simplicity, we may assume that the Sperner's operations have been done to replace all $i$-sets for $i>\frac{n}{2}+2$ and consider $s=\frac{n}{2}+2$. Specifically, at the stage of replacing $(\frac{n}{2}+2)$-sets in Sperner's operations, we chose $\mathscr{A}^{\downarrow}\subseteq (\mathscr{A}-\mathscr{A}^{(\frac{n}{2}+2)})\cup \Delta \mathscr{A}^{(\frac{n}{2}+2)}$ such that $|\mathscr{A}^\diamond|=|\mathscr{A}|$. For ease of argument, let $\mathscr{A}^{\downarrow}= (\mathscr{A}-\mathscr{A}^{(\frac{n}{2}+2)})\cup S$, where $S\subset \Delta \mathscr{A}^{(\frac{n}{2}+2)}$ and $|S|=|\mathscr{A}^{(\frac{n}{2}+2)}|$. That is, $S$ is the set of $(\frac{n}{2}+1)$-sets in $\Delta\mathscr{A}^{(\frac{n}{2}+2)}$ selected as replacements of $\mathscr{A}^{(\frac{n}{2}+2)}$ in Sperner's operations, and $\mathscr{A}_2=\mathscr{A}^{(\frac{n}{2}+1)}\cup S$ (see Figure \ref{figA8.3.1}).
\begin{center}
\begin{tikzpicture}[thick,scale=1]%
\draw(3,10)node[circle, draw, fill=black!100, inner sep=0pt, minimum width=5pt](1_5){};
\draw(4,10)node[circle, draw, fill=black!100, inner sep=0pt, minimum width=5pt](2_5){};
\draw(5,10)node[circle, draw, fill=black!100, inner sep=0pt, minimum width=5pt](3_5){};
\draw(6,10)node[circle, draw, inner sep=0pt, minimum width=5pt](4_5){};

\draw(1,8)node[circle, draw, fill=black!100, inner sep=0pt, minimum width=5pt](1_4){};
\draw(2,8)node[circle, draw, fill=black!100, inner sep=0pt, minimum width=5pt](2_4){};
\draw(3,8)node[circle, draw, fill=black!100, inner sep=0pt, minimum width=5pt](3_4){};
\draw(4,8)node[circle, draw, inner sep=0pt, minimum width=5pt](4_4){};
\draw(5,8)node[circle, draw, fill=black!100, inner sep=0pt, minimum width=5pt, label={[yshift=-0.7cm] 270:$S$}](5_4){};
\draw(6,8)node[circle, draw, fill=black!100, inner sep=0pt, minimum width=5pt](6_4){};
\draw(8,8)node[circle, draw, fill=black!100, inner sep=0pt, minimum width=5pt](7_4){};

\draw(-1,6)node[circle, draw, fill=black!100, inner sep=0pt, minimum width=5pt](1_3){};
\draw(0,6)node[circle, draw, fill=black!100, inner sep=0pt, minimum width=5pt](2_3){};
\draw(1,6)node[circle, draw, inner sep=0pt, minimum width=5pt](3_3){};
\draw(2,6)node[circle, draw, inner sep=0pt, minimum width=5pt](4_3){};
\draw(3,6)node[circle, draw, inner sep=0pt, minimum width=5pt](5_3){};
\draw(4,6)node[circle, draw, fill=black!100, inner sep=0pt, minimum width=5pt](6_3){};
\draw(5,6)node[circle, draw, fill=black!100, inner sep=0pt, minimum width=5pt](7_3){};
\draw(6,6)node[circle, draw, fill=black!100, inner sep=0pt, minimum width=5pt](8_3){};
\draw(7,6)node[circle, draw, fill=black!100, inner sep=0pt, minimum width=5pt](9_3){};
\draw(8,6)node[circle, draw, fill=black!100, inner sep=0pt, minimum width=5pt](10_3){};
\draw(10,6)node[circle, draw, fill=black!100, inner sep=0pt, minimum width=5pt](11_3){};

\draw(1,4)node[circle, draw, inner sep=0pt, minimum width=5pt](1_2){};
\draw(2,4)node[circle, draw, fill=black!100, inner sep=0pt, minimum width=5pt](2_2){};
\draw(3,4)node[circle, draw, fill=black!100, inner sep=0pt, minimum width=5pt](3_2){};
\draw(4,4)node[circle, draw, fill=black!100, inner sep=0pt, minimum width=5pt](4_2){};
\draw(5,4)node[circle, draw, fill=black!100, inner sep=0pt, minimum width=5pt](5_2){};
\draw(6,4)node[circle, draw, fill=black!100, inner sep=0pt, minimum width=5pt](6_2){};
\draw(8,4)node[circle, draw, fill=black!100, inner sep=0pt, minimum width=5pt](7_2){};

\draw(-2,10)node[label={180:$(\frac{n}{2}+2)$-sets}](R5){};
\draw(-2,8)node[label={180:$(\frac{n}{2}+1)$-sets}](R4){};
\draw(-2,6)node[label={180:$(\frac{n}{2})$-sets}](R3){};
\draw(-2,4)node[label={180:$(\frac{n}{2}-1)$-sets}](R2){};

\draw[loosely dotted, line width=0.3mm, >=latex, shorten <= 0.2cm, shorten >= 0.15cm](6.5,8)--(7.5,8);
\draw[loosely dotted, line width=0.3mm, >=latex, shorten <= 0.2cm, shorten >= 0.15cm](8.5,6)--(9.5,6);
\draw[loosely dotted, line width=0.3mm, >=latex, shorten <= 0.2cm, shorten >= 0.15cm](6.5,4)--(7.5,4);

\node[draw, inner xsep=2.5mm,inner ysep=2.5mm, fit=(2_3)(3_3)(4_3)(5_3), label={90:$\mathscr{A}_1$}](A1){};
\node[draw, inner xsep=2.5mm,inner ysep=2.5mm, fit=(1_4)(2_4)(3_4), label={90:$\nabla\mathscr{A}_1$}](shade_A1){};
\node[draw, inner xsep=2.5mm,inner ysep=2.5mm, fit=(5_4)(7_4), label={90:$\Delta\mathscr{A}^{(\frac{n}{2}+2)}$}](shadow_A2){};

\draw[dashed, ->, line width=0.5mm, >=latex, shorten <= 0.2cm, shorten >= 0.15cm, red](1_2)--(1_3);
\draw[dashed, ->, line width=0.5mm, >=latex, shorten <= 0.2cm, shorten >= 0.15cm, red](1_2)--(2_3);
\draw[dashed, ->, line width=0.5mm, >=latex, shorten <= 0.2cm, shorten >= 0.15cm, blue](2_3)--(1_4);
\draw[dashed, ->, line width=0.5mm, >=latex, shorten <= 0.2cm, shorten >= 0.15cm, blue](5_3)--(3_4);

\draw[dotted, ->, line width=0.5mm, >=latex, shorten <= 0.2cm, shorten >= 0.15cm, black!20!green](4_5)--(5_4);
\draw[dotted, ->, line width=0.5mm, >=latex, shorten <= 0.2cm, shorten >= 0.15cm, black!20!green](4_5)--(7_4);

\draw[decorate,decoration={brace, mirror, amplitude=5pt}] (5.6,7.5) -- (8.25,7.5);
\draw(7,7.5)node[label={[yshift=-0.05cm] 270:$\Delta\mathscr{A}^{(\frac{n}{2}+2)}-S$}]{};

\draw[decorate,decoration={brace, mirror, amplitude=5pt}] (4.7,7.5) -- (5.4,7.5);

\draw(7.5,2.5)node[circle, draw, inner sep=0pt, minimum width=5pt, label={[xshift=0.55cm] 0:Elements of $\mathscr{A}$}]{};
\draw[dashed, ->, line width=0.5mm, >=latex, shorten <= 0.2cm, shorten >= 0.15cm](7,2.1)--(8,2.1);
\draw[dashed, ->, line width=0.5mm, >=latex, shorten <= 0.2cm, shorten >= 0.15cm](7,1.9)--(8,1.9);
\draw[dotted, ->, line width=0.5mm, >=latex, shorten <= 0.2cm, shorten >= 0.15cm](7,1.6)--(8,1.6);
\draw[dotted, ->, line width=0.5mm, >=latex, shorten <= 0.2cm, shorten >= 0.15cm](7,1.4)--(8,1.4);
\draw(8,2)node[label={0: Shade}]{};
\draw(8,1.5)node[label={0: Shadow}]{};

\draw(7.25,2.5)node(L1){};
\draw(11,1.5)node(L2){};
\node[draw, inner xsep=2.5mm,inner ysep=2.5mm, fit=(L1)(L2), label={[xshift=-1.5cm] 90:Legend}](){};
\end{tikzpicture}
\captionsetup{justification=centering}
{\captionof{figure}{Sketch of Hasse diagram of $\mathbb{N}_n$; \\to show $\mathscr{A}^{(\frac{n}{2}+2)}=\emptyset$} \label{figA8.3.1}}
\end{center}
\indent\par Since $|\Delta \mathscr{A}^{(\frac{n}{2}+2)}|\ge (1+\frac{6}{n-2})|\mathscr{A}^{(\frac{n}{2}+2)}|>|S|$ by (\ref{eqA5.1.2}), it follows that $\Delta \mathscr{A}^{(\frac{n}{2}+2)}-S\neq\emptyset$. Furthermore, $\mathscr{A}$ is an antichain implies $(\Delta \mathscr{A}^{(\frac{n}{2}+2)}-S)\cap(\nabla\mathscr{A}_1\cup \mathscr{A}^{(\frac{n}{2}+1)})=\emptyset$. Consequently, $\Delta \mathscr{A}^{(\frac{n}{2}+2)}-S$, $\mathscr{A}_2$ and $\nabla\mathscr{A}_1$ are pairwise disjoint families of $(\frac{n}{2}+1)$-sets. It follows that $(\Delta \mathscr{A}^{(\frac{n}{2}+2)}-S)\cup\mathscr{A}_2\cup \nabla\mathscr{A}_1\subseteq {{\mathbb{N}_n}\choose{(n/2)+1}}$. Hence, $|\mathscr{A}_2|<{{n}\choose{(n/2)+1}}-|\nabla\mathscr{A}_1|$, a contradiction to (\ref{eqA5.4.4}). A similar argument holds for $\mathscr{B}$.
\indent\par Next, we show $b(\mathscr{X})\ge\frac{n}{2}$ for $\mathscr{X}=\mathscr{A},\mathscr{B}$. WLOG, suppose $|\mathscr{A}_1|={{n}\choose{n/2}}$. Then, $\mathscr{A}={{\mathbb{N}_n}\choose{n/2}}$ by Sperner's Theorem. It follows that $\mathscr{B}^{(i)}=\emptyset$ for all $i<\frac{n}{2}$. Otherwise, there exist two sets in $\mathscr{A}$ that are each disjoint with any $i$-set in $\mathscr{B}$ with $i<\frac{n}{2}$. By (\ref{eqA5.4.6}), $\mathscr{B}={{\mathbb{N}_n}\choose{(n/2)+1}}$ and we are done.
\indent\par Hence, we may assume $0<|\mathscr{X}_1|<{{n}\choose{n/2}}$ for all $\mathscr{X}=\mathscr{A},\mathscr{B}$. Suppose $\mathscr{A}^{(s)}\neq\emptyset$ for some $s<\frac{n}{2}$. For simplicity, we may assume that the Sperner's operations have been done to replace all $i$-sets for $i>\frac{n}{2}-2$ and consider $s=\frac{n}{2}-1$.
\\
\\Case 1. $\mathscr{A}^{(\frac{n}{2}-1)}={{\mathbb{N}_n}\choose{(n/2)-1}}$. 
\indent\par $\mathscr{A}$ is a antichain implies $\mathscr{A}=\mathscr{A}^{(\frac{n}{2}-1)}$. Then, (\ref{eqA5.4.6}) and Sperner's Theorem imply $\mathscr{B}={{\mathbb{N}_n}\choose{n/2}}$ (and $\kappa^*(k)=0$). Now, every element $X$ of $\mathscr{A}$ is disjoint with ${{|\bar{X}|}\choose{n/2}}={{(n/2)+1}\choose{n/2}}=\frac{n}{2}+1\ge 3$ elements of $\mathscr{B}$, a contradiction.
\\
\\Case 2. $\emptyset\neq\mathscr{A}^{(\frac{n}{2}-1)}\subset{{\mathbb{N}_n}\choose{(n/2)-1}}$. 
\indent\par Recall that at the stage of replacing $(\frac{n}{2}-1)$-sets in Sperner's operations, we chose $\mathscr{A}^{\uparrow}\subseteq(\mathscr{A}-\mathscr{A}^{(\frac{n}{2}-1)})\cup \nabla \mathscr{A}^{(\frac{n}{2}-1)}$ such that $|\mathscr{A}^{\uparrow}|=|\mathscr{A}|$. For ease of argument, let $\mathscr{A}^{\uparrow}= (\mathscr{A}-\mathscr{A}^{(\frac{n}{2}-1)})\cup S$, where $S\subset \nabla \mathscr{A}^{(\frac{n}{2}-1)}$ and $|S|=|\mathscr{A}^{(\frac{n}{2}-1)}|<|\nabla \mathscr{A}^{(\frac{n}{2}-1)}|$; the strict inequality follows from (\ref{eqA5.1.1}). In other words, $\mathscr{A}_1=\mathscr{A}^{(\frac{n}{2})}\cup S$.
\\
\\Subcase 2.1. $\mathscr{A}^{(\frac{n}{2})}=\emptyset$.
\indent\par Note that $\nabla\mathscr{A}_1\subseteq\nabla(\nabla\mathscr{A}^{(\frac{n}{2}-1)})$. Suppose $\nabla\mathscr{A}_1\subset\nabla(\nabla\mathscr{A}^{(\frac{n}{2}-1)})$ (see Figure \ref{figA8.3.2}). Since $\mathscr{A}$ is a antichain, $\nabla\mathscr{A}_1$, $\nabla(\nabla\mathscr{A}^{(\frac{n}{2}-1)})-\nabla\mathscr{A}_1$, and $\mathscr{A}_2$ are disjoint families of $(\frac{n}{2}+1)$-sets. It follows that $\nabla\mathscr{A}_1\cup (\nabla(\nabla\mathscr{A}^{(\frac{n}{2}-1)})-\nabla\mathscr{A}_1) \cup \mathscr{A}_2\subseteq{{\mathbb{N}_n}\choose{(n/2)+1}}$. Hence, $|\mathscr{A}_2|<{{n}\choose{(n/2)+1}}-|\nabla\mathscr{A}_1|$, a contradiction to (\ref{eqA5.4.4}).
\\\\
\begin{center}
\begin{tikzpicture}[thick,scale=1]%
\draw(1,8)node[circle, draw, fill=black!100, inner sep=0pt, minimum width=5pt](1_4){};
\draw(2,8)node[circle, draw, fill=black!100, inner sep=0pt, minimum width=5pt](2_4){};
\draw(3,8)node[circle, draw, fill=black!100, inner sep=0pt, minimum width=5pt](3_4){};
\draw(4,8)node[circle, draw, inner sep=0pt, minimum width=5pt](4_4){};
\draw(5,8)node[circle, draw, inner sep=0pt, minimum width=5pt](5_4){};
\draw(6,8)node[circle, draw, inner sep=0pt, minimum width=5pt](6_4){};
\draw(8,8)node[circle, draw, inner sep=0pt, minimum width=5pt](7_4){};

\draw(-1,6)node[circle, draw, fill=black!100, inner sep=0pt, minimum width=5pt](1_3){};
\draw(0,6)node[circle, draw, fill=black!100, inner sep=0pt, minimum width=5pt](2_3){};
\draw(1,6)node[circle, draw, fill=black!100, inner sep=0pt, minimum width=5pt](3_3){};
\draw(2,6)node[circle, draw, fill=black!100, inner sep=0pt, minimum width=5pt](4_3){};
\draw(3,6)node[circle, draw, fill=black!100, inner sep=0pt, minimum width=5pt](5_3){};
\draw(4,6)node[circle, draw, fill=black!100, inner sep=0pt, minimum width=5pt](6_3){};
\draw(5,6)node[circle, draw, fill=black!100, inner sep=0pt, minimum width=5pt](7_3){};
\draw(6,6)node[circle, draw, fill=black!100, inner sep=0pt, minimum width=5pt](8_3){};
\draw(7,6)node[circle, draw, fill=black!100, inner sep=0pt, minimum width=5pt](9_3){};
\draw(8,6)node[circle, draw, fill=black!100, inner sep=0pt, minimum width=5pt](10_3){};
\draw(10,6)node[circle, draw, fill=black!100, inner sep=0pt, minimum width=5pt](11_3){};

\draw(1,4)node[circle, draw, inner sep=0pt, minimum width=5pt](1_2){};
\draw(2,4)node[circle, draw, inner sep=0pt, minimum width=5pt](2_2){};
\draw(3,4)node[circle, draw, inner sep=0pt, minimum width=5pt](3_2){};
\draw(4,4)node[circle, draw, fill=black!100, inner sep=0pt, minimum width=5pt](4_2){};
\draw(5,4)node[circle, draw, fill=black!100, inner sep=0pt, minimum width=5pt](5_2){};
\draw(6,4)node[circle, draw, fill=black!100, inner sep=0pt, minimum width=5pt](6_2){};
\draw(8,4)node[circle, draw, fill=black!100, inner sep=0pt, minimum width=5pt](7_2){};

\draw(-2,8)node[label={180:$(\frac{n}{2}+1)$-sets}](R4){};
\draw(-2,6)node[label={180:$(\frac{n}{2})$-sets}](R3){};
\draw(-2,4)node[label={180:$(\frac{n}{2}-1)$-sets}](R2){};

\draw[loosely dotted, line width=0.3mm, >=latex, shorten <= 0.2cm, shorten >= 0.15cm](6.5,8)--(7.5,8);
\draw[loosely dotted, line width=0.3mm, >=latex, shorten <= 0.2cm, shorten >= 0.15cm](8.5,6)--(9.5,6);
\draw[loosely dotted, line width=0.3mm, >=latex, shorten <= 0.2cm, shorten >= 0.15cm](6.5,4)--(7.5,4);

\node[draw, inner xsep=2.5mm,inner ysep=2.5mm, fit=(1_3)(2_3)(3_3)(4_3)(5_3), label={90:$\nabla\mathscr{A}^{(\frac{n}{2}-1)}$}](A1){};
\node[draw, inner xsep=2.5mm,inner ysep=2.5mm, fit=(1_4)(2_4)(3_4), label={90:$\nabla(\nabla\mathscr{A}^{(\frac{n}{2}-1)})$}](shade_A1){};
\node[draw, inner xsep=2.5mm,inner ysep=2.5mm, fit=(4_4)(7_4), label={90:$\mathscr{A}_2$}](shade_A1){};

\draw[dashed, ->, line width=0.5mm, >=latex, shorten <= 0.2cm, shorten >= 0.15cm, red](1_2)--(1_3);
\draw[dashed, ->, line width=0.5mm, >=latex, shorten <= 0.2cm, shorten >= 0.15cm, red](3_2)--(5_3);
\draw[dashed, ->, line width=0.5mm, >=latex, shorten <= 0.2cm, shorten >= 0.15cm, blue](1_3)--(1_4);
\draw[dashed, ->, line width=0.5mm, >=latex, shorten <= 0.2cm, shorten >= 0.15cm, blue](5_3)--(3_4);

\draw[decorate,decoration={brace, mirror, amplitude=5pt}] (1.6,7.5) -- (3.25,7.5);
\draw(1.7,7.1)node[label={[yshift=-0.05cm] 0: $\nabla\mathscr{A}_1$}]{};

\draw[decorate,decoration={brace, mirror, amplitude=5pt}] (0.6,5.5) -- (3.25,5.5);
\draw(0.9,5)node[label={[yshift=-0.05cm] 0: $\mathscr{A}_1=S$}]{};

\draw(7.5,2.5)node[circle, draw, inner sep=0pt, minimum width=5pt, label={[xshift=0.55cm] 0:Elements of $\mathscr{A}$}]{};
\draw[dashed, ->, line width=0.5mm, >=latex, shorten <= 0.2cm, shorten >= 0.15cm](7,2.1)--(8,2.1);
\draw[dashed, ->, line width=0.5mm, >=latex, shorten <= 0.2cm, shorten >= 0.15cm](7,1.9)--(8,1.9);
\draw(8,2)node[label={0: Shade}]{};

\draw(7.25,2.5)node(L1){};
\draw(11,2)node(L2){};
\node[draw, inner xsep=2.5mm,inner ysep=2.5mm, fit=(L1)(L2), label={[xshift=-1.5cm] 90:Legend}](){};
\end{tikzpicture}
\captionsetup{justification=centering}
{\captionof{figure}{Sketch of Hasse diagram of $\mathbb{N}_n$;\\Subcase 2.1 with $\nabla\mathscr{A}_1\subset\nabla(\nabla\mathscr{A}^{(\frac{n}{2}-1)})$.} \label{figA8.3.2}}
\end{center}
\indent\par Since $\mathscr{A}$ is an antichain,
\begin{align}
|\mathscr{A}_1|-|\nabla\mathscr{A}_1|=\kappa(|\mathscr{A}_1|)=\kappa^*(k)+\kappa({{n}\choose{(n/2)}})-\kappa^*(k+{{n}\choose{n/2}}-|\mathscr{A}_1|)\ge 0, \label{eqA5.4.7}
\end{align}
where we used (\ref{eqA5.4.5}),  (\ref{eqA5.4.3}), and Lemma \ref{lemA8.3.8} and $\kappa^*(k)\ge 0$ respectively.
\indent\par Now, suppose $\nabla\mathscr{A}_1=\nabla(\nabla\mathscr{A}^{(\frac{n}{2}-1)})$. Recall that the replacements in Sperner's operations are one-to-one and therefore $|\mathscr{A}^{(\frac{n}{2}-1)}|=|\mathscr{A}_1|$. So, by (\ref{eqA5.1.1}) and $\emptyset\neq\mathscr{A}^{(\frac{n}{2}-1)}\subset{{\mathbb{N}_n}\choose{(n/2)-1}})$\begin{align*}
|\nabla\mathscr{A}_1|=|\nabla(\nabla\mathscr{A}^{(\frac{n}{2}-1)})|\ge\frac{n}{n+2} |\nabla\mathscr{A}^{(\frac{n}{2}-1)}|>\frac{n}{n+2}\cdot\frac{n+2}{n}\cdot|\mathscr{A}^{(\frac{n}{2}-1)}| =|\mathscr{A}^{(\frac{n}{2}-1)}|
=|\mathscr{A}_1|,
\end{align*}
which contradicts (\ref{eqA5.4.7}).
\\
\\Subcase 2.2. $\mathscr{A}^{(\frac{n}{2})}\neq\emptyset$.
\\
\\Claim 1: $\nabla\mathscr{A}_1=\nabla(\mathscr{A}^{(\frac{n}{2})}\cup \nabla\mathscr{A}^{(\frac{n}{2}-1)})$.
\indent\par Since $\mathscr{A}_1=\mathscr{A}^{(\frac{n}{2})}\cup S\subset \mathscr{A}^{(\frac{n}{2})}\cup \nabla\mathscr{A}^{(\frac{n}{2}-1)}$, we have  $\nabla\mathscr{A}_1\subseteq \nabla(\mathscr{A}^{(\frac{n}{2})}\cup \nabla\mathscr{A}^{(\frac{n}{2}-1)})$. Now, if $|\nabla\mathscr{A}_1|<|\nabla(\mathscr{A}^{(\frac{n}{2})}\cup \nabla\mathscr{A}^{(\frac{n}{2}-1)})|$, then $|\mathscr{A}_2|\le |{{\mathbb{N}_n}\choose{(n/2)+1}}|-|\nabla(\mathscr{A}^{(\frac{n}{2})}\cup \nabla\mathscr{A}^{(\frac{n}{2}-1)})|<{{n}\choose{(n/2)+1}}-|\nabla\mathscr{A}_1|=|\mathscr{A}_2|$, a contradiction to (\ref{eqA5.4.4}). Hence, the claim follows.
\\
\indent\par Since $\mathscr{A}^{(\frac{n}{2})}\cup S=\mathscr{A}_1$, it follows from Claim 1 that $\nabla(\nabla\mathscr{A}^{(\frac{n}{2}-1)}-S)\subseteq\nabla(\mathscr{A}^{(\frac{n}{2})}\cup S)$. Note that $|\nabla\mathscr{A}^{(\frac{n}{2}-1)}|>|\mathscr{A}^{(\frac{n}{2}-1)}|=|S|$ by (\ref{eqA5.1.1}). Let $A_0\in \nabla\mathscr{A}^{(\frac{n}{2}-1)}-S$, where $A_{i^*}\subset A_0$ for some $A_{i^*}\in \mathscr{A}^{(\frac{n}{2}-1)}$ (see Figure \ref{figA8.3.3}).
\begin{center}
\begin{tikzpicture}[thick,scale=1]%
\draw(1,8)node[circle, draw, fill=black!100, inner sep=0pt, minimum width=5pt](1_4){};
\draw(2,8)node[circle, draw, fill=black!100, inner sep=0pt, minimum width=5pt](2_4){};
\draw(3,8)node[circle, draw, fill=black!100, inner sep=0pt, minimum width=5pt](3_4){};
\draw(4,8)node[circle, draw, inner sep=0pt, minimum width=5pt](4_4){};
\draw(5,8)node[circle, draw, inner sep=0pt, minimum width=5pt](5_4){};
\draw(6,8)node[circle, draw, inner sep=0pt, minimum width=5pt](6_4){};
\draw(8,8)node[circle, draw, inner sep=0pt, minimum width=5pt](7_4){};

\draw(-1,6)node[circle, draw, fill=black!100, inner sep=0pt, minimum width=5pt](1_3){};
\draw(0,6)node[circle, draw, fill=black!100, inner sep=0pt, minimum width=5pt, label={[xshift=-0.15cm] 270:$A_0$}](2_3){};
\draw(1,6)node[circle, draw, fill=black!100, inner sep=0pt, minimum width=5pt, label={[xshift=-0.15cm, yshift=0.1cm] 270:$A^{\uparrow}_{i^*}$}](3_3){};
\draw(2,6)node[circle, draw, fill=black!100, inner sep=0pt, minimum width=5pt](4_3){};
\draw(3,6)node[circle, draw, fill=black!100, inner sep=0pt, minimum width=5pt](5_3){};
\draw(4,6)node[circle, draw, inner sep=0pt, minimum width=5pt](6_3){};
\draw(5,6)node[circle, draw, inner sep=0pt, minimum width=5pt](7_3){};
\draw(6,6)node[circle, draw, fill=black!100, inner sep=0pt, minimum width=5pt](8_3){};
\draw(7,6)node[circle, draw, fill=black!100, inner sep=0pt, minimum width=5pt](9_3){};
\draw(8,6)node[circle, draw, fill=black!100, inner sep=0pt, minimum width=5pt](10_3){};
\draw(10,6)node[circle, draw, fill=black!100, inner sep=0pt, minimum width=5pt](11_3){};

\draw(1,4)node[circle, draw, inner sep=0pt, minimum width=5pt](1_2){};
\draw(2,4)node[circle, draw, inner sep=0pt, minimum width=5pt, label={270:$A_{i^*}$}](2_2){};
\draw(3,4)node[circle, draw, inner sep=0pt, minimum width=5pt](3_2){};
\draw(4,4)node[circle, draw, fill=black!100, inner sep=0pt, minimum width=5pt](4_2){};
\draw(5,4)node[circle, draw, fill=black!100, inner sep=0pt, minimum width=5pt](5_2){};
\draw(6,4)node[circle, draw, fill=black!100, inner sep=0pt, minimum width=5pt](6_2){};
\draw(8,4)node[circle, draw, fill=black!100, inner sep=0pt, minimum width=5pt](7_2){};

\draw(-2,8)node[label={180:$(\frac{n}{2}+1)$-sets}](R4){};
\draw(-2,6)node[label={180:$(\frac{n}{2})$-sets}](R3){};
\draw(-2,4)node[label={180:$(\frac{n}{2}-1)$-sets}](R2){};

\draw[->, line width=0.3mm, >=latex, shorten <= 0.1cm, shorten >= 0.2cm](2_2)--(2_3);
\draw[->, line width=0.3mm, >=latex, shorten <= 0.1cm, shorten >= 0.2cm](2_2)--(3_3);

\draw[loosely dotted, line width=0.3mm, >=latex, shorten <= 0.2cm, shorten >= 0.15cm](6.5,8)--(7.5,8);
\draw[loosely dotted, line width=0.3mm, >=latex, shorten <= 0.2cm, shorten >= 0.15cm](8.5,6)--(9.5,6);
\draw[loosely dotted, line width=0.3mm, >=latex, shorten <= 0.2cm, shorten >= 0.15cm](6.5,4)--(7.5,4);

\draw(2,5.65)node(X0){};
\node[draw, inner xsep=2.5mm,inner ysep=2.5mm, fit=(1_3)(2_3)(3_3)(4_3)(5_3)(X0), label={90:$\nabla\mathscr{A}^{(\frac{n}{2}-1)}$}](A1){};
\node[draw, inner xsep=2.5mm,inner ysep=2.5mm, fit=(1_4)(2_4)(3_4), label={90:$\nabla(\nabla\mathscr{A}^{(\frac{n}{2}-1)})$}](shade_A1){};

\node[draw, inner xsep=2.5mm,inner ysep=2.5mm, fit=(4_4)(7_4), label={90:$\mathscr{A}_2$}](shade_A1){};

\draw[dashed, ->, line width=0.5mm, >=latex, shorten <= 0.2cm, shorten >= 0.15cm, red](1_2)--(1_3);
\draw[dashed, ->, line width=0.5mm, >=latex, shorten <= 0.2cm, shorten >= 0.15cm, red](3_2)--(5_3);
\draw[dashed, ->, line width=0.5mm, >=latex, shorten <= 0.2cm, shorten >= 0.15cm, blue](1_3)--(1_4);
\draw[dashed, ->, line width=0.5mm, >=latex, shorten <= 0.2cm, shorten >= 0.15cm, blue](5_3)--(3_4);
\draw[dashed, ->, line width=0.5mm, >=latex, shorten <= 0.2cm, shorten >= 0.15cm, black!20!green](3_3)--(1_4);
\draw[dashed, ->, line width=0.5mm, >=latex, shorten <= 0.2cm, shorten >= 0.15cm, black!20!green](7_3)--(3_4);
\draw[decorate,decoration={brace, mirror, amplitude=5pt}] (1.6,7.5) -- (3.25,7.5);
\draw(1.7,7.1)node[label={[yshift=-0.05cm] 0: $\nabla\mathscr{A}_1$}]{};

\draw[decorate,decoration={brace, mirror, amplitude=5pt}] (0.7,5.2) -- (3.25,5.2);
\draw(1.45,4.7)node[label={[yshift=-0.05cm] 0: $S$}]{};

\draw[decorate,decoration={brace, mirror, amplitude=10pt}] (0.7,5.1) -- (5.6,5.1);
\draw(2.75,4.5)node[label={[yshift=-0.05cm] 0: $\mathscr{A}_1$}]{};

\draw(7.5,2.5)node[circle, draw, inner sep=0pt, minimum width=5pt, label={[xshift=0.55cm] 0:Elements of $\mathscr{A}$}]{};
\draw[dashed, ->, line width=0.5mm, >=latex, shorten <= 0.2cm, shorten >= 0.15cm](7,2.1)--(8,2.1);
\draw[dashed, ->, line width=0.5mm, >=latex, shorten <= 0.2cm, shorten >= 0.15cm](7,1.9)--(8,1.9);
\draw(8,2)node[label={0: Shade}]{};
\draw[->, line width=0.3mm, >=latex, shorten <= 0.2cm, shorten >= 0.2cm](7,1.5)--(8,1.5);
\draw(8,1.5)node[label={0: Proper subset}]{};

\draw(7.25,2.5)node(L1){};
\draw(11,1.5)node(L2){};
\node[draw, inner xsep=2.5mm,inner ysep=2.5mm, fit=(L1)(L2), label={[xshift=-1.5cm] 90:Legend}](){};
\end{tikzpicture}
\captionsetup{justification=centering}
{\captionof{figure}{Sketch of Hasse diagram of $\mathbb{N}_n$; Subcase 2.2.} \label{figA8.3.3}}
\end{center}
Claim 2: $|T|=k$, (recall $T:=\{i|A^{\uparrow}_i\cap B^{\uparrow}_i=\emptyset\}$), i.e. there are exactly $k$ disjoint pairs $(A^{\uparrow}_i,B^{\uparrow}_i)$, $i=1,2,\ldots,k$, where $A^{\uparrow}_i\in\mathscr{A}_1$ and $B^{\uparrow}_i\in\mathscr{B}_1$.
\indent\par Suppose $|T|<k$. Observe that $\nabla(\nabla\mathscr{A}^{(\frac{n}{2}-1)}-S)\subseteq\nabla(\mathscr{A}^{(\frac{n}{2})}\cup S)$ implies $\mathscr{A}_1\cup \mathscr{A}_2\cup\{A_0\}$ is a antichain. Furthermore, $|A_0|=\frac{n}{2}$ and $\frac{n}{2}\le |B|\le \frac{n}{2}+1$ for all $B\in \mathscr{B}_1\cup\mathscr{B}_2$ implies $A_0\cap B\neq \emptyset$ except for at most one element (if $B=\bar{A_0}\in \mathscr{B}_1$) of $\mathscr{B}_1$. In other words, $\mathscr{A}_1\cup \mathscr{A}_2\cup\{A_0\}$ and $\mathscr{B}_1\cup \mathscr{B}_2$ are cross-intersecting antichains with at most $k$ disjoint pairs $(A^{\uparrow}_i,B^{\uparrow}_i)$ and size more than $|\mathscr{A}|+|\mathscr{B}|$, a contradiction to the optimality of $\mathscr{A}$ and $\mathscr{B}$. Hence, this claim follows.
\\
\indent\par Since $|\bar{A}_{i^*}|=\frac{n}{2}+1$ and by (\ref{eqA5.4.4}), either $\bar{A}_{i^*}\in\mathscr{B}_2$ or $\bar{A}_{i^*}\in\nabla\mathscr{B}_1$. Suppose $\bar{A}_{i^*}\in\mathscr{B}_2$ holds, i.e. $A_{i^*}$ is disjoint with some $(\frac{n}{2}+1)$-set $B\in\mathscr{B}_2$, $\bar{A}_{i^*}=B$. However, the replaced pair $(A^{\uparrow}_{i^*}, B^{\downarrow}_{i^*})$ of this disjoint pair $(A_{i^*},B)$ is intersecting since $|A^{\uparrow}_{i^*}|+|B^{\downarrow}_{i^*}|>n$. This contradicts the fact of having exactly $k$ disjoint pairs before and after Sperner's operations by Claim 2. Hence, $\bar{A}_{i^*}\in\nabla\mathscr{B}_1$. Now, $A_{i^*}$ is disjoint with at most one element of $\mathscr{B}$ (and hence $\mathscr{B}_1$ too) implies $\mathscr{B}_1\cap\Delta \{\bar{A}_{i^*}\}=\{B^{\uparrow}_{i^*}\}$, with $A^{\uparrow}_{i^*}=\bar{B}^{\uparrow}_{i^*}$.
\\
\\Claim 3: $A_0\cap Y\neq\emptyset$ for all $Y\in \mathscr{B}_1\cup\mathscr{B}_2$.
\indent\par Suppose not. Since $|A_0|=\frac{n}{2}$, $A_0$ is disjoint with some $B^{\uparrow}_j\in\mathscr{B}_1$ for some $j$. Furthermore, $A_0\neq A^{\uparrow}_{i^*}=\bar{B}^{\uparrow}_{i^*}$ implies $j\neq i^*$. That is, $A_{i^*}\in \mathscr{A}$ is originally disjoint with two distinct sets in $\mathscr{B}$, namely $B_{i^*}$ and $B_j$, a contradiction. Hence, this claim follows.
\\
\indent\par By Claim 3, $\mathscr{A}_1\cup \mathscr{A}_2\cup\{A_0\}$ and $\mathscr{B}_1\cup \mathscr{B}_2$ are cross-intersecting antichains with $k$ disjoint pairs and size more than $|\mathscr{A}|+|\mathscr{B}|$, a contradiction to the maximality of $\mathscr{A}$ and $\mathscr{B}$.
\indent\par A similar argument holds for $\mathscr{B}$.
\begin{flushright}
$\Box$
\end{flushright}

\section{Exactly $k$ disjoint pairs}
\indent\par We end this paper by considering the variation closely related to Theorem \ref{thmA5.2.5}. Let $n\ge 4$ be an even integer and $\mathscr{A}$ and $\mathscr{B}$ be two antichains of $\mathbb{N}_n$. Suppose for some integer $k\le \min\{|\mathscr{A}|, |\mathscr{B}|\}$, and for all $A_i\in \mathscr{A}$, $B_j\in\mathscr{B}$, $A_i\cap B_j=\emptyset$ if and only if $i=j\le k$, i.e. there are \textit{exactly} $k$ disjoint pairs. Determine the maximum $|\mathscr{A}|+|\mathscr{B}|$.
\indent\par As mentioned, Corollary \ref{corA5.2.7} provides a tight upper bound if $\kappa(k)=\kappa^*(k)$. We are interested to know if (\ref{eqA5.2.1}) holds in the case of $\kappa(k)<\kappa^*(k)$. We think this is true but do not have a proof. 
\indent\par We remark that we may choose $A^{\uparrow}_i$ and $B^{\uparrow}_i$ so that $|T|=k$ in the proof of Thoerem \ref{thmA5.2.5}. Indeed, if $i\neq j$, then $A^{\uparrow}_i\neq A^{\uparrow}_j$ (and $B^{\uparrow}_i\neq B^{\uparrow}_j$ resp.). Otherwise, $A^{\uparrow}_i=A^{\uparrow}_j=\bar{B}^{\uparrow}_j$ implies $A_i\cap B_j=\emptyset$, a contradiction.
\indent\par Hence, following the proof outline of Theorem \ref{thmA5.2.5}, it seems that the Conjecture \ref{conA5.5.1}, which is analogous to our auxiliary Proposition \ref{ppnA5.2.4}, will be useful. We verified the conjecture for small even integers of $n$ but a proof remains elusive.
\begin{con}\label{conA5.5.1}
Let $n\ge 4$ be an even integer, and $a$ and $k$ be integers such that $0\le a,k \le {{n}\choose{n/2}}$. Then, 
\begin{align*}
\kappa_{\frac{n}{2}}({{n}\choose{n/2}})+\kappa_{\frac{n}{2}}(k)\le \kappa_{\frac{n}{2}}(a)+\kappa_{\frac{n}{2}}^*(k+{{n}\choose{n/2}}-a).
\end{align*}
\end{con}
\textbf{Acknowledgement}
\indent\par The first author would like to thank the National Institute of Education, Nanyang Technological University of Singapore, for the generous support of the Nanyang Technological University Research Scholarship.

\end{document}